\documentclass[a4paper,11pt,leqno]{article}
\usepackage{amsmath,amsfonts,amsthm,amssymb}
\usepackage[left=3cm, right=3cm, top=3cm, bottom=3.5cm]{geometry}
\usepackage{graphicx}
\usepackage{mathrsfs}
\usepackage{tikz}
\numberwithin{equation}{section}
\usepackage[final]{pdfpages}

\theoremstyle{plain}
\newtheorem{theorem}{Theorem}[section]
\newtheorem{lemma}[theorem]{Lemma}

\theoremstyle{definition}

\newcommand{\R}{\mathbb{R}}
\newcommand{\N}{\mathbb{N}}

\newcommand{\supp}{\operatorname{supp}}

\newcommand{\cA}{\mathscr{A}}
\newcommand{\cN}{\mathscr{N}}
\newcommand{\cV}{\mathscr{V}}
\newcommand{\cH}{\mathscr{H}}
\usepackage{color}

\title{{\Large On the structure of the infinitesimal generators of scalar \\
one-dimensional semigroups with discrete Lyapunov functionals}}

\author{Giorgio Fusco\\
        Dipartimento di Matematica Pura ed Applicata \\
				Universit\`a degli Studi dell'Aquila \\
        Via Vetoio, 67010 Coppito, L'Aquila, ITALY \\
        {\tt fusco@univaq.it}\\
        {}
  \and
        Carlos Rocha\\
        Instituto Superior T\'ecnico -- Universidade de Lisboa\\
        Avenida Rovisco Pais, 1049--001 Lisboa, PORTUGAL\\
        {\tt crocha@tecnico.ulisboa.pt}\\
        {\tt http://camgsd.tecnico.ulisboa.pt}\\
        {}  }

\date{version of \today}

\begin{document}

\maketitle

\setlength{\parindent}{0cm}

\begin{abstract}

Dynamical systems generated by scalar reaction-diffusion equations on an interval enjoy special properties 
that lead to a very simple structure for the semiflow. Among these properties, the monotone behavior of the 
number of zeros of the solutions plays an essential role. This discrete Lyapunov functional contains 
important information on the spectral behavior of the linearization and leads to a Morse-Smale description 
of the dynamical system.
Other systems, like the linear scalar delay differential equations under monotone feedback conditions, 
possess similar kinds of discrete Lyapunov functions.

Here we discuss and characterize classes of linear equations that generate semiflows acting on $C^0[0,1]$ 
or on  $C^1[0,1]$ which admit discrete Lyapunov functions related to the zero number. 
We show that, if the space is $C^1[0,1]$, the corresponding equations are essentially parabolic partial 
differential equations. In contrast, if the space is $C^0[0,1]$, the corresponding equations are 
generalizations of monotone feedback delay differential equations.

\end{abstract}

\section{Introduction} \label{sec1}

A large amount of the research conducted on the theory of nonlinear dynamical
systems has addressed the study of semiflows possessing certain types
of discrete Lyapunov functionals. This particularity is shared by certain systems
generated either by scalar parabolic partial differential equations on an interval or by scalar
delay differential equations in one space variable.
These systems, or better, their linearizations along the flow, possess a certain
integer valued functional (often called a zero number or a lap number) which is a
nonincreasing function of time. This feature has shown to be essential for many of
the special properties that these systems have in common. Such studies have been
conducted by many authors and, for references to the decreasing property of the zero
number, see \cite{mat82,ang88,mal88,mase96a,mase96b}. See also \cite{lou19} for a more
recent extension of this result.
Regarding the far reaching consequences of this property to the study of the
corresponding infinite dimensional dynamical systems, see
\cite{hen85,ang86,brfi88,brfi89,fima89,firo96,firo00,krwa01}.
Here we also mention the equivalent results for ordinary differential equations
obtained in \cite{fuol88,fuol90} and the result for maps \cite{okm93}.

For references on the general theory of dynamical systems generated by parabolic
partial differential equations or by delay differential equations we refer to
\cite{hal88,hmo02,rob01,seyo02}.

In these notes, we address the problem of characterizing classes of linear semigroups acting on 
$C^0[0,1]$ or on $C^1[0,1]$ (with suitable boundary conditions) which admit discrete Lyapunov 
functionals related to the zero number. We will see that this characterization depends on the 
particular discrete functional considered and  also on the smoothness of the function space.

In the next Section, we recall the decreasing property of the zero number in the setting of linear
parabolic differential equations. Then, we recollect the discrete Lyapunov functions, also derived
from the zero number, for scalar delay differential equations under monotone feedback conditions.

We present, in the third Section, our Theorem \ref{th31} as a main result. It consists of a characterization
of the infinitesimal generator of a $C_0$-semigroup of bounded linear operators in $C^1[0,1]$ with the
zero number decreasing property. We show that, in the appropriate setting, this infinitesimal generator
corresponds necessarily to a second order differential operator.

In Section \ref{sec4} we prepare to also consider the characterization of the infinitesimal generators of
semiflows with different discrete Lyapunov functions $V^\mp$ derived from the zero number like in
the case of delay differential equations. Our main result in this Section, Theorem 4.1, yields 
a characterization of generators of semigroups in $C^1[0,1]$ with such discrete Lyapunov
functionals.  We show that the discrete Lyapunov functions $V^\mp$ are not specific of delay differential
equations. In fact, in Theorem 4.1, we prove that these discrete Lyapunov functions determine,
in the appropriate settings, classes of infinitesimal generators corresponding to second order
differential operators under non-separated boundary conditions.

In Section \ref{sec5} we illustrate this result with a degenerate example showing that linear scalar
delay equations naturally fit this description.

Finally, in Section \ref{sec6} we consider the characterization of generators of semigroups on $C^0[0,1]$
with discrete Lyapunov functions $V^\mp$. Our main result, Theorem \ref{th61}, shows that these semigroups
include delay differential equations. We complete  Theorem \ref{th61} with a Remark where we discuss how the 
result in Theorem \ref{th31} changes when we relax the smoothness of the function space by replacing $C^1[0,1]$ 
with $C^0[0,1]$. 

The notion of zero or lap number has been also extended and used for the description of the dynamics of 
fully non linear scalar parabolic equations \cite{lap18,lap22} and to the case where the differential 
operator is the $p$-Laplacian \cite{gebr07}.

The global dynamics of semiflows that preserve an ordering of the function space is known to have various 
interesting properties  \cite{hir88,pol89,smth91}. For the case of equation \eqref{eq201} where a complete 
description of the set of equilibria and the proof of the Morse-Smale property are possible the key point 
is that the zero number induces a structure of total order in the function space.

\vspace*{.3cm}

\textbf{Acknowledgments.}
We wish to thank the anonymous referees whose excellent comments and suggestions
greatly improved the presentation of this paper.

CR dedicates this paper to Giorgio Fusco, a colleague and coauthor, in appreciation of his longtime friendship.
Both authors wish to express their gratitude to the memory of Professor George Sell which
in 2012 insisted that we should publish these results.
This work was partially supported by FCT/Portugal through projects UID/MAT/04459/2019 and UIDB/04459/2020.

\section{Discrete Lyapunov functions} \label{sec2}

We next recall the remarkable decreasing property of the zero number in the setting of linear
parabolic differential equations. Loosely speaking, this property
can be described by saying that if $t \mapsto u(\cdot,t)$ is a solution of a linear
parabolic differential equation with the corresponding boundary conditions, then the number
$z(u(\cdot,t))$ of zeros of the function $u(\cdot,t)$ is a nonincreasing function of $t$.
A rigorous description is the following Theorem \ref{th21} (see \cite{ang88,lou19}), where
for simplicity we choose Neumann boundary conditions by taking
$C^1_n[0,1] = \left\{ \varphi \in C^1[0,1]: \varphi'(0)=\varphi'(1)=0 \right\}$.

Let
\begin{equation} \label{eq201}
u_t = a(x) u_{xx} + b(t,x) u_x + c(t,x) u  \ , \ x \in (0,1) \ .
\end{equation}
If the coefficients $a, b, c$ are sufficiently smooth and $a$ is positive, problem
\eqref{eq201} generates an evolution operator $S(t,t_0):C^1_n[0,1]\rightarrow C^1_n[0,1], t\ge t_0$,
corresponding to the solutions $S(t,t_0) u_0(x) = u(t,x)$ with $u(t_0,x)=u_0(x)$.
If, in addition, $b,c$ are independent of $t$ then \eqref{eq201} generates a linear $C_0$-semigroup
$\{T(t)=S(t,0)\}_{t\ge 0}$. See, for example, \cite{paz83,hen81}.

For $\varphi\in C^1_n[0,1]$ let the {\em zero number} $z(\varphi)\in\N_0\cup\{\infty\}$ denote the number
of strict sign changes of $x\mapsto\varphi(x)$. The precise definition, appearing in \cite{ang88} for
$\varphi\in C^0[0,1]$, see also \cite{mase96a}, is the following. Let
\begin{equation} \label{eq202}
\begin{split}
z(\varphi) = & \sup \{ k\geq 0: \exists\;\,\{x_i\}_{i=0}^k\subset[0,1]\;\text{such that} \; x_{i-1}<x_i\; \\
& \text{and} \; \varphi( x_{i-1})\varphi(x_i)<0\;,\; i=1,\ldots,k \} \ .
\end{split}
\end{equation}
If $\varphi\in C^0[0,1]$ satisfies $\varphi(x)\geq 0$ or $\varphi(x)\leq 0$
for $x\in[0,1]$ the zero number is set as $z(\varphi)=0$.
Of course this definition also holds for $\varphi\in C^1_n[0,1]$.
Notice that, according to this definition, $z(\varphi) = 0$ for $\varphi\equiv 0$.

Then, $z(S(t,t_0)\varphi)$ is a monotone nonincreasing function of $t>t_0$, \cite{ang88,lou19}
and see also \cite{firo96}. More specifically let
$\cN=\left\{ \varphi \in C^1_n[0,1]): \varphi(x)=0 \Rightarrow \varphi'(x) \ne 0 \right\}$ be the open
dense subset of functions with all zeros nondegenerate. Restricted to $\cN$ the zero number $z$ is
the map which associates to each $\varphi \in \cN$ the number of zeros of $\varphi$. Note that
$\cN = \cup_{k=0}^\infty \cN_k$, $\cN_k = \left\{ \varphi \in \cN : z(\varphi)=k \right\}$ and that $z|_\cN$
is continuous and locally constant. Then

\begin{theorem} \label{th21}
For $\varphi \in C^1_n[0,1]$, $\varphi \not\equiv 0$,
\begin{enumerate}
    \item [(i)] the set $\Theta = \left\{t \in (t_0,+\infty): S(t,t_0)\varphi\not \in \cN \right\}$ is a finite set;
    \item [(ii)] for $t \in \Theta$ there exists a positive $\varepsilon_0$ such that for all
		$0<\varepsilon<\varepsilon_0$,
    \begin{equation} \label{eq203}
    z(S(t+\varepsilon,t_0)\varphi) < z(S(t-\varepsilon,t_0)\varphi) \ .
    \end{equation}
\end{enumerate}
\end{theorem}

This nonincreasing character of $z$ along solutions of \eqref{eq201} is essential
for the description of the global dynamics of semilinear parabolic equations of the form
\begin{equation} \label{eq204}
u_t = a(x) u_{xx} + f(x,u,u_x)  \ , \ x \in (0,1) \ ,
\end{equation}
with the appropriate boundary conditions. See \cite{brfi88,brfi89,furo91,firo96,firo00} for the
pioneering results. See also \cite{lap22} for the fully nonlinear case.
For a smooth nonlinearity $f\in C^2([0,1]\times\R^2,\R)$ satisfying suitable dissipative conditions, the
semilinear reaction-diffusion equation \eqref{eq204} defines a global semiflow in the space $X=H^1(0,1)$.
Then, $z$ provides a discrete Lyapunov function for the linearized flow around any given solution of 
\eqref{eq204} or for the evolution of the difference between any two solutions $u_1, u_2$, of \eqref{eq204}. 
In fact,
\begin{equation} \label{eq205}
t \mapsto z(u_1(t,\cdot)-u_2(t,\cdot))
\end{equation}
is nonincreasing for $t>0$ since $u_1-u_2$ satisfies a variational equation of the form \eqref{eq201}. 
Moreover, considering the difference $u(t+\tau,\cdot)-u(t,\cdot)$ we obtain that 
\begin{equation} \label{eq205a}
t \mapsto z(u_t(t,\cdot))
\end{equation} 
is also nonincreasing for $t>0$, see \cite{fima89} and Lemma 4.5 in \cite{mana97}.
The existence of a lap number for degenerate parabolic equations is considered in \cite{gebr07}.

Different examples of discrete Lyapunov functions are provided by scalar delay differential equations of the form
\begin{equation} \label{eq206}
\dot x(t) = h(x(t),x(t-1)) \ ,
\end{equation}
with $h\in C^2(\R^2,\R)$ and monotone feedback conditions
\begin{equation} \label{eq207}
h_v(u,v)<0 \mbox{ or } h_v(u,v)>0 \ .
\end{equation}
These equations define semiflows $x_t=T(t)x_0, t\ge 0$, in $X=C^0[-1,0]$ by
$x_t(\theta)=x(t+\theta), \theta\in[-1,0]$.

The following functionals rooted on $z$ are defined on $X$ for each feedback condition \eqref{eq207},
respectively
\begin{equation} \label{eq208}
V^-(\varphi) = 2 \left\lfloor\frac{z(\varphi)}{2}\right\rfloor+1 \quad , \quad
V^+(\varphi) = 2 \left\lfloor\frac{z(\varphi)+1}{2}\right\rfloor \ ,
\end{equation}
where $\lfloor\cdot\rfloor$ denotes the floor function.
Notice that for $\varphi\in X$, $V^-(\varphi)$ is always odd while $V^+(\varphi)$ is always even.

Then, for each feedback case, $V^\mp$ provides a discrete Lyapunov function for the difference
of any two solutions $x^1_t, x^2_t$, since
\begin{equation} \label{eq209}
t \mapsto V^\mp(x^1_t-x^2_t)
\end{equation}
is nonincreasing for $t>0$. For references, see \cite{mase96a} and \cite{nie23}.
These examples will be relevant for the characterization of infinitesimal generators in
our final three Sections.

\section{Generators of semiflows with zero number decay} \label{sec3}

It is natural to ask if the property of possessing a discrete Lyapunov functional $z$
in the sense of Theorem \ref{th21} is characteristic of differential second order operators.
To answer this question, our main result is the following theorem.

\begin{theorem} \label{th31}
Let $A$ be the infinitesimal generator of a $C_0$-semigroup of
bounded linear operators $\{T(t)\}_{t\ge 0}$, $T(t): C^1_n[0,1] \rightarrow C^1_n[0,1]$,
with $D(A)=C^2[0,1]\cap C^1_n[0,1]$ and such that
\begin{enumerate}
    \item [(i)] the set $\Theta = \left\{ t \in (0,+\infty): T(t) \varphi\not\in\cN \right\}$
		is a finite set for every $\varphi \in C^1_n[0,1]$;
		\item [(ii)] for all $0 < t_1 \le t_2$, with $t_1,t_2 \not\in\Theta$ the following holds
    \begin{equation}\label{eq301}
		z(T(t_1)\varphi) \ge z(T(t_2)\varphi) \ .
    \end{equation}
\end{enumerate}
Then there exist $\alpha, \gamma \in C^0[0,1]$ with $\alpha$ nonnegative, and bounded $\beta \in C^0(0,1)$
such that for all $\varphi \in D(A)$ we have
\begin{equation} \label{eq302}
\left(A\varphi\right)(x) = \alpha(x) \varphi_{xx}(x) + \beta(x) \varphi_x(x) +
\gamma(x) \varphi(x) \quad , \quad 0 < x < 1 \ .
\end{equation}
\end{theorem}

Here we assume Neumann boundary conditions, but similar results can be obtained for the cases
of Dirichlet, Robin or even periodic boundary conditions.
In Section 4 we show that this characterization extends also to other types of non-separated
boundary conditions.
\vskip.2cm
The main point in the proof of Theorem \ref{th31}, see Lemma \ref{le33} below, concerns the analysis 
of the behavior of $A$ on functions with degenerate zeros. We prepare with an auxiliary lemma.

\begin{lemma}\label{le32}
Given $\xi\in [0,1]$, assume $\varphi: [0,1] \rightarrow \R$, $\psi: [0,1] \rightarrow \R$, are $C^2$ 
functions such that
\begin{enumerate}
    \item [(i)] $\varphi(\xi)=\varphi'(\xi)=\psi(\xi)=\psi'(\xi)=\psi''(\xi)=0$;
    \item [(ii)] $\varphi''(\xi)>0$;
    \item [(iii)] $\varphi'(x) \ne 0$ for every $x \in (0,1)\setminus\{\xi\}$.
\end{enumerate}
For each $\lambda\in\R$ define $\varphi^\lambda: [0,1]\rightarrow\R$ by
\begin{equation}\label{eq303}
\varphi^\lambda(x) = \varphi(x) + \lambda \psi(x) \ , \ x\in [0,1] \ .
\end{equation}
Then the set of $\lambda\in\R$ such that $\varphi^\lambda$ has a degenerate zero in $[0,1]\setminus\{\xi\}$
has Lebesgue measure zero.
\end{lemma}

\begin{proof}
By (ii) and (iii) it follows that $\varphi(x)>0$ for every $x \ne \xi$.
Assume first that $x\in\{0,1\}, x \ne\xi$. Then there is at most one value of $\lambda$ such that
$\varphi^\lambda(x)=0$.
It remains to show that the set of $\lambda$ such that $\varphi^\lambda$ has a degenerate zero in
$(0,1)\setminus \{\xi\}$ has Lebesgue measure zero.
To show this we note that $x \in (0,1)\setminus \{\xi\}$ is a degenerate zero of $\varphi^\lambda$ 
if and only if
\begin{equation} \label{eq304}
\left\{\begin{array}{l}
\varphi(x)+\lambda\psi(x)=0,\\ \\
\varphi^\prime(x)+\lambda\psi^\prime(x)=0
\end{array}\right.
\;\Leftrightarrow\quad\; \frac{\psi(x)}{\varphi(x)}= \frac{\psi^\prime(x)}{\varphi^\prime(x)}= -\frac{1}{\lambda} \ ,
\end{equation}
where we have used assumption (iii) and the positivity of $\varphi(x)$ for $x\ne\xi$ which, in particular,
imply $\lambda\ne 0$. 

Consider now the function $\Psi:(0,1)\setminus\{\xi\}\rightarrow\R$ defined by
$\Psi(x)=\frac{\psi(x)}{\varphi(x)}$. Since $\varphi$ and $\varphi^\prime$ do not vanish for $x\ne\xi$,
$\Psi$ is well defined, of class $C^1$ and
\begin{equation} \label{eq305}
\Psi^\prime(x)=\frac{\varphi'(x)}{\varphi(x)}\left(\frac{\psi'(x)}{\varphi'(x)}-\Psi(x)\right) \ ,
\end{equation}
vanishes if and only if $\Psi(x)=\frac{\psi'(x)}{\varphi'(x)}$.
From this and \eqref{eq304} we conclude that $x\in(0,1)\setminus\{\xi\}$ is a degenerate zero of $\varphi^\lambda$
if and only if $-\lambda^{-1}$ is a critical value of $\Psi$. Since $\Psi$ is a $C^1$ function defined in an
open set an application of Sard's theorem concludes the proof, see ex.2 pag.170 in \cite{dieu72}.
\end{proof}

\begin{figure}
\begin{center}
\begin{tikzpicture}[xscale=1,yscale=.9]
\draw [] (0,2.5) -- (6.3,2.5);
\draw [] (0,0) -- (6.3,0);

\draw [blue] (1.5,2.5) to [out=0,in=230] (2,2.85);
\draw [blue] (1.5,2.5) to [out=180,in=310] (1,2.85);
\draw [blue] (1,2.85) to [out=130,in=0] (0,3.3);
\draw [blue] (2,2.85) to [out=50,in=180] (6,4.1);

\draw [] (0,-.95) -- (0,1.4);
\draw [](0,2.2)--(0,4);
\draw [](0,-3)--(0,-1.3);

\draw [blue] (1.5,0) to [out=45,in=180] (6,1.5);
\draw [blue] (1.5,0) to [out=225,in=0] (0,-.9);

\draw [] (0,-3) -- (6.3,-3);

\draw [blue] (0,-1.5)--(6,-1.5);

\draw [dotted] (6,2.5)--(6,4.1);
\draw [dotted] (6,0)--(6,1.5);
\draw [dotted] (6,-3)--(6,-1.5);

\node[right] at (6.3,2.5) {$x$};
\node[right] at (6.3,0) {$x$};
\node[right] at (6.3,-3) {$x$};

\node[left] at (.35,2.25) {$0$};
\node[left] at (.35,-.25) {$0$};
\node[left] at (.35,-3.25) {$0$};

\node[right] at (5.8,2.2) {$1$};
\node[right] at (5.8,-.295) {$1$};
\node[right] at (5.8,-3.25) {$1$};

\draw [fill] (1.5,2.5) circle [radius=0.025];
\draw [fill] (1.5,0) circle [radius=0.025];
\draw [fill] (1.5,-3) circle [radius=0.025];

\node [below] at (1.5,2.5) {$\xi$};
\node [below] at (1.5,0) {$\xi$};
\node [below] at (1.5,-3) {$\xi$};


\node [right] at (0,3.8) {$a^\xi$};
\node [right] at (0,1.2) {$b^\xi$};
\node [right] at (0,-1.3) {$c^\xi$};
\end{tikzpicture}
\end{center}
\caption{The functions $a^\xi, b^\xi$ and $c^\xi$.}
\label{fig1}
\end{figure}
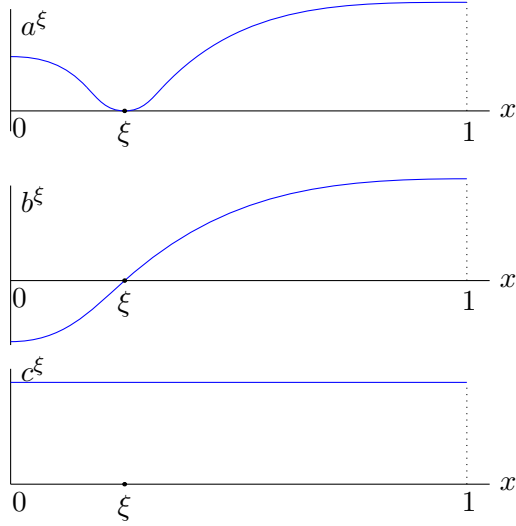

Next we associate to each $\xi\in [0,1]$ three functions
$a^\xi, b^\xi, c^\xi \in D(A)$ such that (here we use the partial derivative notation $v'=v_x$,
see Figure \ref{fig1})
\begin{equation} \label{eq306}
\left\{ \begin{array}{l}
a^\xi(\xi) = a^\xi_x(\xi) = b^\xi(\xi) = 0 \ , \mbox{ for every }\xi\in[0,1] \ , \\  \\
a^\xi_{xx}(\xi) = c^\xi(\xi) = 1 \ , \mbox{ for every }\xi\in[0,1] \ , \\ \\
b^\xi_x(\xi) = 1 \ , \mbox{ for every } \xi \in (0,1) \ , \\ \\
a^\xi_x(x)(x-\xi) > 0 \ , \mbox{ for every } x \in (0,1)\setminus \{\xi\} \ .
\end{array} \right.
\end{equation}

\begin{lemma} \label{le33}
Let $\xi \in [0,1]$ and $\psi^\xi \in D(A)$ be such that
\begin{equation} \label{eq307}
\psi^\xi(\xi) = \psi^\xi_x(\xi) = \psi^\xi_{xx}(\xi) = 0 \ .
\end{equation}
Then, if $A$ is the generator of a semigroup that satisfies the assumptions
in Theorem \ref{th31}, it results that
\begin{equation} \label{eq308}
\left(A\psi^\xi\right)(\xi) = 0 \ .
\end{equation}
\end{lemma}

\begin{proof}
Let $\varphi^{\xi,\lambda} = a^\xi + \lambda \psi^\xi$. Then, on the basis of Lemma \ref{le32}, there is a
sequence $\{\lambda_k\}_{k=1}^\infty$ such that
\begin{equation} \label{eq309}
\lambda_k \lambda_{k+1} < 0 \quad , \quad \lim_{k\rightarrow +\infty} |\lambda_k| = +\infty
\end{equation}
and such that all zeros of $\varphi^{\xi,\lambda_k}$ are nondegenerate
except the one at $\xi$ which is the unique degenerate zero of $\varphi^{\xi,\lambda_k}$.
From the condition on the second derivative of $a^\xi$ in (\ref{eq306})
it also follows that, for each $k=1,2,\dots$, there is a $\delta_k>0$ such
that for every $x \in (\xi-\delta_k,\xi+\delta_k)$ we have
\begin{equation} \label{eq310}
\varphi^{\xi,\lambda_k}(x) \ge \frac{1}{4}(x-\xi)^2 \ .
\end{equation}

This and the fact that the zeros of $\varphi^{\xi,\lambda_k}$ outside $(\xi-\delta_k,\xi+\delta_k)$
are nondegenerate imply that the number of zeros of $\varphi^{\xi,\lambda_k}$ is finite and equal to
$N+1$, where $N$ is the number of zeros of $\varphi^{\xi,\lambda_k}$ in
$[0,1]\setminus (\xi-\delta_k,\xi+\delta_k)$. Therefore
\begin{equation} \label{eq311}
\varphi^{\xi,\lambda_k}\in\partial\cN_N \ .
\end{equation}
Let $\Theta_k$ be the set defined in (i) in Theorem \ref{th31} with $\psi=\varphi^{\xi,\lambda_k}$. 
Since $\Theta_k$ is a finite set, there exists $t_k>0$ such that $(0,t_k]\cap\Theta_k=\emptyset$. 
It follows that, for $t\in(0,t_k]$, $T(t)\varphi^{\xi,\lambda_k}$ has only nondegenerate zeros and, since
 $\lim_{t\rightarrow 0^+}T(t)\varphi^{\xi,\lambda_k}=\varphi^{\xi,\lambda_k}$, has exactly $N$ zeros in 
$[0,1]\setminus(\xi-\delta_k,\xi+\delta_k)$. Hence

\begin{equation} \label{eq312}
\begin{array}{l}
T(t_k)\varphi^{\xi,\lambda_k}\in\cN_{\bar N}, \ \text{ for some } \ \bar N\ge N \,\\ \\
\text{and}\quad T(t_k)\varphi^{\xi,\lambda_k}(\xi)<0,\quad\Rightarrow\quad\;\bar N\geq N+2 \ .
\end{array}
\end{equation}
 
{\em Claim:} Assume $(A\psi^\xi)(\xi)\neq 0$. Then $k$ large and $t_k>0$ sufficiently small imply
\begin{equation}\label{eq313}
\lambda_k(A\psi^\xi)(\xi)(T(t_k)\varphi^{\xi,\lambda_k})(\xi)>0 \ .
\end{equation}

From the assumption that $a^\xi, \psi^\xi \in D(A)$ it follows
\begin{equation} \label{eq314}
\left\{ \begin{array}{l}
T(t)a^\xi - a^\xi - t Aa^\xi = o_1(t) \\ \\
T(t)\psi^\xi - \psi^\xi - t A\psi^\xi = o_2(t)
\end{array} \right. \ ,
\end{equation}
where $o_i(t) \in C^1(0,1), i=1,2$, are small order terms of $t>0$,
that is, $\|\frac{o_i(t)}{t}\|\rightarrow 0$ as $t\searrow 0$.
Hence $t_k>0$ sufficiently small and  $\left(A\psi^\xi\right)(\xi)\ne 0$ yield

\begin{equation} \label{eq315}
\left|\frac{o_1(t_k)}{t_k}(\xi)\right| <1 \ ,\quad \left|\frac{o_2(t_k)}{t_k}(\xi)\right| < 
\frac{1}{2}\vert(A\psi^\xi)(\xi)\vert \ .
\end{equation}

This, $\varphi^{\xi,\lambda_k}(\xi)=0$ and
\begin{equation} \label{eq316}
\begin{array}{cl}
&\left(T(t_k)\varphi^{\xi,\lambda_k}\right)(\xi)=
\left(T(t_k)\varphi^{\xi,\lambda_k}\right) (\xi) -
\varphi^{\xi,\lambda_k} (\xi) \\ \\
&= t_k\left((Aa^\xi)(\xi) +
\lambda_k (A\psi^\xi)(\xi)+ (o_1(t_k)/t_k)(\xi) +
\lambda_k(o_2(t_k)/t_k)(\xi)\right) \ , \end{array}
\end{equation}
imply that, for $k$ sufficiently large,
$\left(T(t_k)\varphi^{\xi,\lambda_k}\right)(\xi)$ has the same sign as
$\lambda_k\left(A\psi^\xi\right)(\xi)$ and the claim is established.

Since $\cN_{\bar N}$ is open we can choose an open neighborhood 
$V_k\subset\cN_{\bar N}$ of $T(t_k)\varphi^{\xi,\lambda_k}$ so small that
\begin{equation} \label{eq317}
\mathrm{sign}(v(\xi))=\mathrm{sign}((T(t_k)\varphi^{\xi,\lambda_k})(\xi))
=\mathrm{sign}(\lambda_k(A\psi^\xi)(\xi)) \ , \ v\in V_k \ ,
\end{equation}
where we have also used \eqref{eq313}.
The continuity of $T(t_k)$ implies $T(t_k)U_k\subset V_k$ for some open nonempty neighborhood $U_k$ of
$\varphi^{\xi,\lambda_k}$ and \eqref{eq311} implies in particular $T(t_k)O_k\subset V_k$ where we have set
$O_k=U_k\cap\cN_N$.
Then \eqref{eq309} and \eqref{eq317} imply that $k$ can be chosen so that
\begin{equation} \label{eq318}
\left(T(t_k)\varphi^{\xi,\lambda_k}\right)(\xi) < 0 \ .
\end{equation}

Then \eqref{eq313} yields
\begin{equation} \label{eq319}
N=z(\varphi)<z(T(t_k)\varphi)=\bar N \ , \ \varphi\in O_k \ ,
\end{equation}
in contradiction with assumption (ii) of Theorem \ref{th31}. This contradiction proves Lemma \ref{le33}.
\end{proof}

Now, based on Lemma \ref{le33}, we complete the proof of Theorem \ref{th31}.
Given $x \in [0,1]$ and $\varphi \in D(A)$ we write 
\begin{equation} \label{eq320}
\varphi = \varphi(x) c^x + \varphi_x(x) b^x + \varphi_{xx}(x) a^x + \psi^{x,\varphi}
\end{equation}
where $a^x,b^x$ and $c^x$ are the functions defined in \eqref{eq306}. Since $\varphi$ and  $a^x,b^x,c^x$ 
belong to $D(A)$ also $\psi^{x,\varphi} \in D(A)$ and \eqref{eq306} implies 
\begin{equation} \label{eq321}
\psi^{x,\varphi}(x) = \psi^{x,\varphi}_x(x) = \psi^{x,\varphi}_{xx}(x) = 0 \ .
\end{equation}

Therefore, Lemma \ref{le33} implies $\left(A\psi^{x,\varphi}\right)(x)=0$ and it follows that
\begin{equation} \label{eq322}
(A\varphi)(x) = \varphi(x) (Ac^x)(x) + \varphi_x(x) (Ab^x)(x) + \varphi_{xx}(x) (Aa^x)(x) \ , \ x \in [0,1] \ ,
\end{equation}
which is \eqref{eq302} with $\alpha, \beta, \gamma$ defined by
\begin{equation} \label{eq323}
\gamma(x)=(Ac^x)(x) \ , \ \beta(x)=(Ab^x)(x) \ , \ \alpha(x)=(Aa^x)(x) \ , \ x\in[0,1] \ .
\end{equation}
On the basis of Lemma \ref{le33}, $(Aa^\xi)(\xi)=\alpha(\xi)<0$ will lead again to  \eqref{eq318}.
This proves that $\alpha$ is a nonnegative function. 
It remains to prove the smoothness of $\alpha, \gamma$. 

\vskip.2cm

Consider the constant map $\varphi \equiv 1$ in $D(A)$. Then, from the abstract theory of $C_0$-semigroups 
on $C_n^1[0,1]$ we have $A1 \in C_n^1[0,1]$. Hence, by \eqref{eq322} and \eqref{eq323} 
\begin{equation} \label{eq324}
\gamma=A 1\in C_n^1[0,1] \ .
\end{equation}

Given $x_0\in(0,1)$ let $\varphi,\psi\in D(A)$ be such that:
\begin{equation} \label{eq325}
\begin{split}
\varphi(x_0)=0 \ , \ \varphi^\prime(x_0)=1 \ , \ \varphi^{\prime\prime}(x_0)=0 \ ,\\
\psi(x_0)=0 \ , \ \psi^\prime(x_0)=0 \ , \ \psi^{\prime\prime}(x_0)=1 \ .
\end{split}
\end{equation}
Set $d(x)=\varphi^\prime(x)\psi^{\prime\prime}(x)-\psi^\prime(x)\varphi^{\prime\prime}(x)$ and observe that 
$d(x_0)=1$. This and $\varphi,\psi\in D(A)$ imply that there is $\delta_0>0$ such that 
\begin{equation}  \label{eq326}
\left.\begin{array}{l}
d(x)\geq\frac1{2} \ , \\ \\
\vert\varphi^\prime-1\vert\leq\frac1{2} \ , \quad\vert\psi^{\prime\prime}-1\vert\leq\frac{1}{2} \ , \\ \\
\vert\psi^\prime\vert\leq\frac1{2} \ , \quad\vert\varphi^{\prime\prime}\vert\leq\frac1{2} \ ,
\end{array}\right. \quad \text{for} \ x\in(x_0-\delta_0,x_0+\delta_0) \ .
\end{equation}
With this choice of $\varphi$ and $\psi$, from (3.22) and (3.23) we get a system of two equation in $\beta(x)$ 
and $\alpha(x)$ that can be solved for $x\in(x_0-\delta_0,x_0+\delta_0)$ yielding:
\begin{equation} \label{327}
\begin{split}
\beta(x)=\frac{((A\varphi)(x)-\varphi(x)\gamma(x))\psi^{\prime\prime}(x)-
((A\psi)(x)-\psi(x)\gamma(x))\varphi^{\prime\prime}(x)}
{\varphi^\prime(x)\psi^{\prime\prime}(x)-\psi^\prime(x)\varphi^{\prime\prime}(x)} \ , \\ \\
\alpha(x)=\frac{((A\psi)(x)-\psi(x)\gamma(x))\varphi^\prime(x)-
((A\varphi)(x)-\varphi(x)\gamma(x))\psi^\prime(x)}
{\varphi^\prime(x)\psi^{\prime\prime}(x)-\psi^\prime(x)\varphi^{\prime\prime}(x)} \ .
\end{split}
\end{equation} 
From this and the continuity of $A\varphi$ and $A\psi$ it follows that $\beta$ and $\alpha$ are continuous 
in $(x_0-\delta_0,x_0+\delta_0)$ and therefore in $(0,1)$ since $x_0\in(0,1)$ is arbitrary. 
From \eqref{eq326} we also have that $\beta$ and $\alpha$ are bounded in $(0,1)$.
It remains to discuss the behavior of $\alpha$ for $x\rightarrow 0,1$. To this end chose $\varphi\in D(A)$ 
such that $\varphi(0)=\varphi^\prime(0)=0$ and $\varphi^{\prime\prime}(0)=1$. 
Then, for $x$ in a neighborhood of $0$ it results
\begin{equation} \label{328}
(A\varphi)(x)=\varphi(x)\gamma(x)+\varphi^\prime(x)\beta(x)+\varphi^{\prime\prime}(x)\alpha(x) \ . 
\end{equation} 
This and the fact that $\alpha$ and $\beta$ are bounded implies the existence of the limit
\begin{equation} \label{329}
\lim_{x\rightarrow 0^+}\alpha(x)=(A\varphi)(0) \ , 
\end{equation}
and we can continuously extend $\alpha$ to the interval $[0,1)$ by setting $\alpha(0)=(A\varphi)(0)$. 
Similarly we obtain an extension of $\alpha$ to the interval $[0,1]$ by setting $\alpha(1)=(A\varphi)(1)$.

This concludes the proof of Theorem \ref{th31}. $\hfill \square$

\vskip.2cm

We can ask if the result $\alpha\geq 0$ in Theorem \ref{th31} can be upgraded to $\alpha$ strictly positive. 
Our guess is that this is not possible under only the assumptions of Theorem \ref{th31}.
See Section \ref{sec5} and Remark 6.3 at the end of Section \ref{sec6}. 

\vskip.2cm 

As we have seen, the zero number is well defined in $C^0[0,1]$. An interesting issue is to understand 
the role of the regularity of the phase space on the structure of semigroups and of the corresponding 
generators that enjoy the decay property of a discrete Lyapunov functional.
Since, by lowering the regularity, the class of motions on which the decay property must hold is enlarged, 
we expect that changing $C_n^1[0,1]$ to $C^0[0,1]$ restricts the class of semigroups that admit $z$ as 
a discrete Lyapunov functional. For a discussion of this question we refer to Remark 6.2.

\section{Generators of semigroups in $C^1[0,1]$ with discrete Lyapunov functions $V^\mp$} \label{sec4}

In this Section we consider the problem of characterizing the generators of $C_0$-semigroups $T(t)$, 
$t\geq 0$, which admit the discrete Lyapunov functional $V^-$ or $V^+$ defined in \eqref{eq208}.
We assume that the function space is $C^1[0,1]$ endowed with suitable boundary conditions.  We observe 
that the arguments developed in the proofs of Lemma \ref{le32} and  Lemma \ref{le33} apply also to the case 
where the discrete Lyapunov functional is $V^-$ or $V^+$. 
This is a consequence of the fact that, if \eqref{eq308} is violated, as shown in Lemma \ref{le33}, 
there are $t_k$ and $\lambda_k$ such that \eqref{eq312} holds and the number of zeros increases of at 
least by two contradicting the decay property of $z$ and $V^-$ and $V^+$. 
From \eqref{eq308} it follows that again the generator is a second order operator and the problem is to 
understand what restrictions to the boundary conditions derive from the assumption that the semigroup admits 
$V^-$ or $V^+$ as a discrete Lyapunov functional, see Theorem \ref{th41} below.

\vskip.2cm

In Section \ref{sec2} we have introduced the functionals $V^\pm$ in connection with the scalar delay 
equation \eqref{eq206} in the phase space $C^0[-1,0]$ which is the standard choice for equation \eqref{eq206}. 
Obviously, with the appropriate definition of the map $\theta\rightarrow x_t(\theta)$ we can replace 
$C^0[-1,0]$ with the space $C^0[0,1]$. The fact is that the right space for the delay equation is a 
space of continuous functions. 
In contrast, as stated in Theorem \ref{th41}, the generator of a semigroup with the same discrete 
Lyapunov functional, acting on a function space of $C^1$ functions, is a second order operator. 
This suggests that the smoothness of the function space plays an important role on the characterization 
of generators with discrete Lyapunov functionals. We confirm this in Section \ref{sec6} where we show 
that changing $C^1$ to $C^0$ implies that the generator has $\alpha\equiv 0$ and corresponds to a 
generalized kind of delay equation.

We consider the space
\begin{equation} \label{eq401}
C^1_{ns}[0,1] = \left\{ \varphi \in C^1[0,1]: B_0(\varphi)=B_1(\varphi)=0 \right\} \ ,
\end{equation}
where $B_0,B_1:C^1[0,1]\rightarrow\R$ are the boundary operators
\begin{equation} \label{eq402}
\left\{ \begin{array}{l}
B_0(\varphi)= \varphi'(0) + \delta_{00}\varphi(0) + \delta_{01}\varphi(1) \ , \\ \\
B_1(\varphi)= \varphi'(1) + \delta_{10}\varphi(0) + \delta_{11}\varphi(1) \ .
\end{array} \right.
\end{equation}
For this type of boundary value problems we refer to ch.11,12 of \cite{cole72}.
See also \cite{lou19} for results regarding the zeros at the boundary.
We remark that $\delta_{01}=\delta_{10}=0$ correspond to the separated Robin boundary conditions not
under consideration here.

Our result regarding the characterization of the infinitesimal generators is the following theorem.

\begin{theorem} \label{th41}
Let $A$ be the infinitesimal generator of a $C_0$-semigroup of bounded linear operators 
$\{T(t)\}_{t\ge 0}, T(t): C^1_{ns}[0,1] \rightarrow C^1_{ns}[0,1]$ with
$D(A) = C^2[0,1] \cap C^1_{ns}[0,1]$ and such that $T(t)$ admits the discrete Lyapunov functional
defined by $V^-$ (or $V^+$, alternatively) in \eqref{eq207}.

Then, there exist $\alpha, \gamma \in C^0[0,1]$, with $\alpha$ nonnegative, and bounded $\beta\in C^0(0,1)$  
such that for all $\varphi \in D(A)$ we have
\begin{equation} \label{eq403}
\left(A\varphi\right)(x) = \alpha(x) \varphi_{xx}(x) + \beta(x) \varphi_x(x) +
\gamma(x) \varphi(x) \ , \ 0 < x < 1 \ .
\end{equation}
Furthermore, the cross-boundary constants in \eqref{eq402} satisfy
\begin{equation} \label{eq404}
\delta_{01} < 0 \ , \ \delta_{10} > 0 \qquad (\mbox{or } \delta_{01} > 0 \ , \
\delta_{10} < 0, \mbox{ respectively}) \ .
\end{equation}
\end{theorem}

\begin{proof}
We focus on the case of $V^-$. The discussion of the case $V^+$ is similar. We let $\cN\subset C^1_{ns}[0,1]$ 
denote the open dense subset of functions with no zeros on the boundary of $[0,1]$ and all (interior) zeros 
nondegenerate.
We remark that in this case we have $\cN = \cup_{k=0}^\infty \cN_{2k+1}$, with
$\cN_{2k+1}=\left\{ \varphi\in\cN : V^-(z(\varphi))=2k+1 \right\}$.

We first consider the case of degenerate interior zeros and address later the zeros on the boundary.

As we have already observed  Lemma \ref{le32} holds as stated and we can consider the family 
$\varphi^{\xi,\lambda}=\varphi^\xi+\lambda\psi^\xi$ with $\varphi^\xi$ and $\psi^\xi$ as in Lemma \ref{le33}. 
It results
\begin{equation} \label{eq405}
\varphi^{\xi,\lambda} \in \partial\cN_{2k+1} \ ,
\end{equation}
for some $k$ and, proceeding as in the proof of Lemma \ref{le33}, we have that $(A\psi^\xi)(\xi)\neq 0$ 
implies the existence of $t>0$ and $\lambda\in\R$ such that $(T(t)\varphi^{\xi,\lambda})(\xi)<0$. 
It follows that $z((T(t)\varphi^{\xi,\lambda})(\xi))=2k^\prime+1\geq 2k+1+2$ and, as in the proof of 
Lemma \ref{le33}, this implies the existence of open sets $O\in\cN_{2k+1}$ and $V\in\cN_{2k'+1}$
such that, for all $\varphi\in O$
\begin{equation} \label{eq406}
2k+1 = V^-(z(\varphi)) < V^-(z(T(t)\varphi)) = 2k'+1 \ ,
\end{equation}
with $T(t)\varphi\in V$, contradicting the decreasing property of the Lyapunov function of the semigroup $T(t)$.
This shows that  $(A\psi^\xi)(\xi)=0$ and the form \eqref{eq403} of the infinitesimal generator $A$ 
follows as in the proof of Theorem \ref{th31}.

To prove \eqref{eq404} we now deal with the zeros on the boundary.
Then, let $\varphi\in C^1_{ns}[0,1]$ denote a function with only nondegenerate zeros and which is zero
at $x=0$ or $x=1$. Note that, due to the boundary conditions \eqref{eq401},\eqref{eq402}, nondegeneracy implies
that these zeros cannot occur simultaneously on both boundaries.

The cross-boundary conditions \eqref{eq404} essentially prevent zeros to occur on the boundary when
$z(\varphi)$ is even. In fact, in this case, the cross-boundary values of $\varphi'$ and $\varphi$ have
the wrong sign, since \eqref{eq401},\eqref{eq402}, imply
\begin{equation} \label{eq407}
\varphi(0)=0 \ \Rightarrow \ \varphi'(0)=-\delta_{01}\varphi(1) \ ,
\end{equation}
\begin{equation} \label{eq408}
\varphi(1)=0 \ \Rightarrow \ \varphi'(1)=-\delta_{10}\varphi(0) \ .
\end{equation}
Then, if $z(\varphi)$ is even, $\varphi(0)=0$ forces $\varphi(1)$ and $\varphi'(0)$ to have
opposite signs, contrary to \eqref{eq407}, \eqref{eq404}.
Likewise, $\varphi(1)=0$ forces $\varphi(0)$ and $\varphi'(1)$ to have the same sign,
contrary to \eqref{eq408}, \eqref{eq404}.
Hence, zeros on the boundaries can only occur when $z(\varphi)$ is odd, in which case $V^-(T(t)\varphi)$ is
constant in some small interval $t\in[0,\varepsilon)$ and $T(t)\varphi\in\cN$ for small $t>0$. 

The proof of smoothness of $\alpha, \beta, \gamma$ follows the same lines as in the proof of Theorem \ref{th31}.
This completes the proof.
\end{proof}

\section{Parabolic and delay operators as members of the same class determined by the decay of $V^-$.} \label{sec5}

Theorem \ref{th41} describes a quite large class of operators parameterized by the functions 
$\alpha\geq 0$, $\beta, \gamma\in\R$, and for the discrete Lyapunov functional $V^-$, 
with the constants $\delta_{01}<0$ and $\delta_{10}>0$.

It is interesting to observe that this class includes second order operators as well as linear delay equations. 
Indeed the operator $A$ corresponding to $\alpha, \gamma\equiv 0$, $\beta\equiv 1$ and $\delta_{10}=p>0$, 
$\delta_{11}=q\in\R$ is given by 
\begin{equation}\label{eq501}
\begin{split}
& (A\phi)(x)=\phi^\prime(x) \ ,\ x\in[0,1) \ , \\ \\
& (A\phi)(1)=-p\phi(0)+q\phi(1) \ .
\end{split}
\end{equation}
This is the infinitesimal generator of the $C_0$-semigroup defined by the linear delay equation
\begin{equation}\label{eq502}
\dot{y}(t)=-py(t-1)+qy(t) \ ,
\end{equation}
see Section 7.1 in \cite{hal77}. As shown in \cite{mase96a} \eqref{eq502} admits the discrete Lyapunov 
functional $V^-$.

\vskip.2cm

We also observe that on the basis of these considerations we cannot hope to deduce the strict 
inequality $\alpha>0$ simply from the decay property of $V^-$.
It is appropriate here to question the role of the cross-boundary conditions \eqref{eq404} and compare with
the motivating examples provided by the scalar delay differential equations under monotone feedback conditions.
The corresponding semiflows also preserve the decreasing properties for the Lyapunov functions defined by $V^\mp$. 
Our illustration example indicates that the cross-boundary conditions in these equations, e.g. $p>0$, 
play the role of the monotone feedback conditions in delay equations.

\section{Generators of semigroups on $C^0[0,1]$ with discrete Lyapunov functionals $V^\mp$} \label{sec6}

Let $T(t): C^0[0,1]\rightarrow C^0[0,1]$, $t\geq 0$, be a $C_0$-semigroup and let $A$ be the infinitesimal 
generator of the semigroup with $D(A)\subset C^0[0,1]$.
We present some results on the structure imposed on $A$ by the assumption that the semigroup admits the
discrete Lyapunov functional $V^-$ or $V^+$.

For $k\in\{0,1\}$ we set $C_{00}^k=\left\{ \varphi\in C^k[0,1]:\supp(\varphi)\subset(0,1) \right\}$, 
and prove the following

\begin{theorem} \label{th61}
Assume that
\begin{enumerate}
\item  $T(t), t\geq 0$ is a $C_0$-semigroup acting on $C^0[0,1]$.
\item $V^\mp(T(t)\varphi)\leq V^\mp(T(\tau)\varphi), \ t\ge\tau, \ \varphi\in C^0[0,1]$.
\item The domain $D(A)\subset C^0[0,1]$ of the infinitesimal generator $A$ of $T(t)$ is a subset of
$C^1[0,1]$ and contains $C_{00}^1$.
\item If $\{\varphi_k\}_{k=1}^\infty\subset D(A)$ is a sequence that converges to $\varphi$ in $C^1[0,1]$, then
$\varphi\in D(A)$ and
\begin{equation} \label{eq601}
\lim_{k\rightarrow+\infty}\|A\varphi_k-A\varphi\|_{C^0[0,1]}=0 \ .
\end{equation}
\end{enumerate}
Then, there exist functions $b,c\in C^0[0,1]$ such that
\begin{equation} \label{eq602}
(A\varphi)(x)=b(x)\varphi^\prime(x)+c(x)\varphi(x) \ , \ x\in[0,1] \ , \ \varphi\in D(A).
\end{equation}
Moreover, if $b(x)\neq 0$ for $x\in[0,1]$, then there is $\widetilde\alpha\in\R$ and $\mp a>0$, for $V^\mp$ 
respectively, such that
\begin{equation} \label{eq603}
\begin{split}
& b>0\Rightarrow D(A)= \left\{ \varphi\in C^1[0,1]:a\varphi(0)+\widetilde\alpha\varphi(1) = 
b(1)\varphi^\prime(1)+c(1)\varphi(1) \right\}, \\ \\
& b<0\Rightarrow D(A)= \left\{ \varphi\in C^1[0,1]:a\varphi(1)+\widetilde\alpha\varphi(0) = 
b(0)\varphi^\prime(0)+c(0)\varphi(0) \right\}.
\end{split}
\end{equation}
\end{theorem}

In the proof we consider the case where (ii) holds with the negative sign. The analysis of the other case is
similar.  We divide the proof in various lemmas.

\begin{lemma} \label{le62}
The assumptions in Theorem \ref{th61} imply:
\begin{equation} \label{eq604}
\supp(A\varphi)\subset\supp(\varphi), \ \varphi\in D(A) \ .
\end{equation}
\end{lemma}

\begin{proof}
We first prove \eqref{eq604} under the assumption $z(\varphi)<+\infty$.
Assume instead that there exists $\varphi\in D(A)$ and an open interval $I=(\alpha,\beta)$ such that
\begin{equation} \label{eq605}
\bar{I}\cap\supp(\varphi)=\emptyset \ \text{and} \ (A\varphi)(x)\neq 0, \ x\in\bar{I} \ .
\end{equation}

Let $I_i=(\alpha_i,\beta_i)\subset I$, $i=1,2,3$ open intervals such that $\beta_1<\alpha_2$ and 
$\beta_2<\alpha_3$.

\begin{figure}
  \begin{center}
\begin{tikzpicture}[scale=1]
\draw[very thin] (-5,0)--(5,0);
\draw[blue, domain=0.4:1] plot (\x, {-.15*(1-cos(600*(\x-.4))});
\draw[blue, domain=-1:-.4] plot (\x, {-.15*(1-cos(600*(\x+1))});
\draw[blue] (-.4,0)--(.4,0);
\draw[blue] (1,0)--(1.5,0);
\draw[blue] (-1,0)--(-1.5,0);
\draw[blue] (4.5,0)--(5,0);
\draw[blue] (-4.5,0)--(-5,0);
\draw[blue, domain=1.5:3] plot (\x, {-.6*(1-cos(120*(\x-1.5))});
\draw[blue, domain=3:5] plot (\x, {-1.2+2*(1-cos(15*(\x-3))});
\draw[dotted] (-5,0)--(-5,-.5);
\draw[dotted] (5,0)--(5,-1);
\draw[blue, domain=-3.75:-1.5] plot (\x, {1.2*(1-cos(120*(\x+4.5))});
\draw[blue, domain=-5:-3.75] plot (\x, {1.2*(1+3*(\x+3.75))+2*(1-cos(90*(\x+3.75))});
\draw[->] (0,-.7)--(.61,-.285);
\draw[->] (0,-.7)--(-.61,-.285);
\path [fill=white] (0,-.7) circle [radius=.16];;
\node[] at (0,-.7) {$\psi$};
\node[below] at (-1.25,0) {$\alpha$};
\node[below] at (1.25,0) {$\beta$};
\draw[] (1.25,-.1)--(1.25,.1);
\draw[] (-1.25,-.1)--(-1.25,.1);
\draw[] (-1.25,.2)--(1.25,.2);
\draw[dotted](-1.8,.2)--(-1.25,.2);
\draw[dotted](1.8,.2)--(1.25,.2);
\node[above] at (0,.4) {$A\varphi$};
\draw[->] (0,.4)--(0,.23);
\draw[->] (2,-.805)--(2.25,-.6);
\node[] at (1.8,-.805) {$\varphi$};
\draw[->] (-1.8,1.4)--(-2.2,1);
\node[] at  (-1.6,1.4) {$\varphi$};
\node[above] at  (-5,0) {$0$};
\node[above] at  (5,0) {$1$};
\end{tikzpicture}
\end{center}
\caption{The maps $\varphi$, $\psi$ and $A\varphi$.}
\label{fig2}
\end{figure}
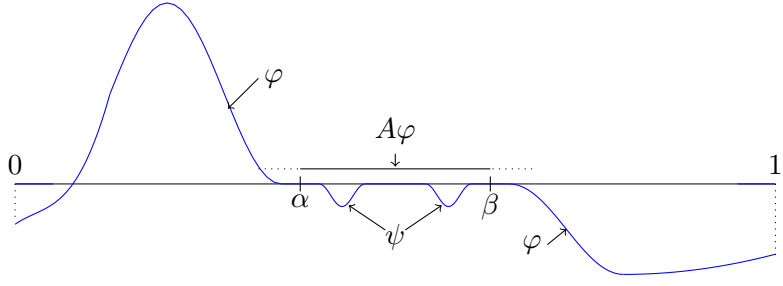

Choose a map $\psi\in C_{00}^1$ with $\supp(\psi)=\bar{I}_1\cup\bar{I}_3$ and such that (see Figure \ref{fig2})
\begin{equation} \label{eq606}
\psi(x)(A\varphi)(x)<0, \ x\in I_1\cup I_3 \ .
\end{equation}
For $\lambda\neq 0$ define $\varphi_\lambda=\lambda\varphi+\psi$ and observe that
\begin{equation} \label{eq607}
\begin{split}
& z(\varphi_1)\leq z(\varphi)+2<+\infty \ , \\
& z(\varphi_\lambda)=z(\varphi_1) \ , \ \lambda\neq 0 \ .
\end{split}
\end{equation}
Fix $\bar{\lambda}>0$ so that
\begin{equation} \label{eq608}
\bar{\lambda}\min_{x\in\bar{I}_2}\vert(A\varphi)(x)\vert>\max_{x\in\bar{I}_2}\vert(A\psi)(x)\vert \ .
\end{equation}
Then we have
\begin{equation} \label{eq609}
\lim_{t\rightarrow 0^+}(T(t)\varphi_{\bar{\lambda}})(x)=\psi(x) \ , \ x\in\bar{I}_1\cup\bar{I}_3 \ ,
\end{equation}
and, for small $t>0$,
\begin{equation} \label{eq610}
\begin{split}
(T(t)\varphi_{\bar{\lambda}})(x)(T(t)\varphi)(x)&=t^2A\varphi_{\bar{\lambda}}A\varphi+\emph{o}(t^2) \\
&=t^2(\bar{\lambda}(A\varphi)^2+A\psi A\varphi)+\emph{o}(t^2)>0 \ , \ x\in\bar{I}_2 \ .
\end{split}
\end{equation}
For $0<t<<1$ this implies that, when $x$ describes the interval $I$, $(T(t)\varphi_{\bar{\lambda}})(x)$
changes sign twice and we have
\begin{equation} \label{eq611}
z(T(t)\varphi_{\bar{\lambda}})=z(\varphi_{\bar{\lambda}})+2 \ ,
\end{equation}
which, in contradiction with (ii), yields $V^-(T(t)\varphi_{\bar{\lambda}})>V^-(\varphi_{\bar{\lambda}})$.
This contradiction concludes the proof for the case $z(\varphi)<+\infty$.

\vskip.2cm

For the general case we construct a sequence $\{\varphi_k\}_{k=1}^\infty\subset D(A)$ which converges to
$\varphi$ in $C^1[0,1]$ and satisfies
\begin{equation} \label{eq612}
\begin{split}
& \supp\varphi_k\subset\supp\varphi \ , \\
& z(\varphi_k)<+\infty \ .
\end{split}
\end{equation}
Assume that there exists a non empty open interval $I$ such that $\bar{I}\cap\supp\varphi=\emptyset$ and
$\vert(A\varphi)(x)\vert> 0$, for $x\in\bar{I}$. From assumption (iv) in Theorem \ref{th61}, we have that
for $k$ sufficiently large $\vert(A\varphi_k)(x)\vert>0$, for $x\in\bar{I}$.
This is in contradiction with the first part of the proof which on the basis of \eqref{eq612} implies
$\supp(A\varphi_k)\subset\supp(\varphi)$. This contradiction shows that \eqref{eq604} holds true also
if $z(\varphi)=+\infty$. It remains to construct the sequence  $\{\varphi_k\}_{k=1}^\infty$.

Let $E=\left\{ x\in[0,1]:\varphi(x)=0 \right\}$ be the set of the zeros of $\varphi$ and $\widetilde{E}$ 
the set of the accumulation points of $E$. For $x\in\widetilde{E}$ we have
$\varphi(x)=\varphi^\prime(x)=0$ and therefore, for $\delta>0$ sufficiently small, we obtain
\begin{equation} \label{eq613}
\left.\begin{array}{l}
\vert\varphi(y)\vert=\emph{o}(\delta) \ , \\ \\
\vert\varphi^\prime(y)\vert\leq f(\delta) \ ,
\end{array} \right. \ 0\leq\vert y-x\vert\leq 2\delta \ , \ y\in[0,1] \ ,
\end{equation}
where $s\rightarrow f(s)$ is a positive function that converges to zero as $s\rightarrow 0^+$.
From \eqref{eq613} it follows that, given $\epsilon>0$, for each $x\in\widetilde{E}$ there exists
$\delta_{x,\epsilon}>0$ such that
\begin{equation} \label{eq614}
\left.\begin{array}{l}
\vert\varphi(y)\vert\leq\emph{o}(\delta_{x,\epsilon})\leq\epsilon\delta_{x,\epsilon} \ , \\ \\
\vert\varphi^\prime(y)\vert\leq\epsilon \ ,
\end{array} \right. \quad 0\leq\vert y-x\vert\leq 2\delta_{x,\epsilon} \ , \ y\in[0,1] \ .
\end{equation}

Since $\widetilde{E}$ is compact and $\{(x-\delta_{x,\epsilon},x+\delta_{x,\epsilon})\cap[0,1]\}_{x\in\widetilde{E}}$
is a open covering of $\widetilde{E}$, there exist an integer $N$ and $x_j\in\widetilde{E}$, $j=1,\ldots,N$ such that
$\widetilde{E}\subset\cup_{j=1}^N[0,1]\cap(x_j-\delta_{x_j,\epsilon},x_j+\delta_{x_j,\epsilon})$.
Set $I_{i,j}=[0,1]\cap(x_j-i\delta_{x_j,\epsilon},x_j+i\delta_{x_j,\epsilon})$, $i=1,2$, $j=1,\ldots,N$ and let
$a_{i,j},b_{i,j}$ be the extremes of $I_{i,j}$.

To conclude the proof of the existence of the sequence $\{\varphi_k\}_{k=1}^\infty$ we define a map
$\varphi_\epsilon\in D(A)\subset C^1[0,1]$ with only a finite number of zeros and such that
$\|\varphi_\epsilon-\varphi\|_{C^1[0,1]}\leq C\epsilon$ with $C>0$ independent of $\epsilon$. We define 
$\varphi_\epsilon$ in each connected component of the set $\cup_{j=1}^NI_{2,j}$. 
Let $I$ be one of these connected components and let $p,q\in\{1,\ldots,N\}$ be defined by 
$a_{1,p}=\min_{x_j\in I}a_{1,j}$ and  $b_{1,q}=\max_{x_j\in I}b_{1,j}$.
Define $\bar\alpha,\bar\beta\in I$ by setting
\begin{equation} \label{eq615}
\begin{split}
&\bar\alpha=\left\{
\begin{array}{ll}
a_{1,p} \ , \ & \text{if} \ a_{2,p}> 0 \ , \\ \\
\min_{x\in\widetilde{E}\cup[a_{1,p},x_p]} \ , & \text{if} \ a_{2,p}=0 \ ,
\end{array}\right. \\ \\
&\bar\beta=\left\{
\begin{array}{ll}
b_{1,q} \ , \ & \text{if} \ b_{2,q}< 1 \ , \\ \\
\max_{x\in\widetilde{E}\cup[x_q,b_{1,q}]} \ , & \text{if} \ b_{2,q}=1 \ .
\end{array}\right.
\end{split}
\end{equation}
 Set
\begin{equation} \label{eq616}
\varphi_\epsilon(x)=0 \ , \ x\in[\alpha,\beta] \ .
\end{equation}
To complete the definition of $\varphi_\epsilon$ in $I$ we focus on the cases $a_{2,p}>0$ and $a_{2,p}=0$, the
discussion of the cases $b_{2,q}<1$ and $b_{2,q}=1$ is analogous. We observe that $a_{2,p}>0$ implies
$a_{2,p}=x_p-2\delta_{x_p,\epsilon}$, $\bar\alpha=x_p-\delta_{x_p,\epsilon}$. Then, define $\varphi_\epsilon$ in the
interval $[x_p-2\delta_{x_p,\epsilon},x_p-\delta_{x_p,\epsilon})$ by setting
\begin{equation} \label{eq617}
\begin{array}{ll}
\varphi_\epsilon(x)=\varphi(x) \ , \ & x\in[0,x_p-2\delta_{x_p,\epsilon}) \ , \\ \\
\varphi_\epsilon(x)=\chi\left((\delta_{x_p,\epsilon})^{-1}(x-x_p+2\delta_{x_p,\epsilon})\right)\varphi(x) \ , \
& x\in[x_p-2\delta_{x_p,\epsilon},x_p-\delta_{x_p,\epsilon}) \ ,
\end{array}
\end{equation}
where $\chi:\R\rightarrow\R$ is a $C^\infty$ nonnegative map that satisfies $\chi(s)=1$, $s\in(-\infty,0)$ and
$\chi(s)=0$, $s\in[1,+\infty]$. From \eqref{eq616} and \eqref{eq617} it follows that to estimate the $C^1$ norm of
$\varphi_\epsilon-\varphi$ in $I\setminus(\bar\beta,1]$ it suffices to look at the interval
$[x_\eta-2\delta_{x_\eta,\epsilon},x_\eta-\delta_{x_\eta,\epsilon})$ where
\begin{equation} \label{eq618}
\varphi_\epsilon(x)-\varphi(x) =
\left(\chi({\delta_{x_\eta,\epsilon}}^{-1}(x-x_\eta+2\delta_{x_\eta,\epsilon}))-1\right) \varphi(x) \ .
\end{equation}

From this expression and \eqref{eq613} we obtain
\begin{equation} \label{eq619}
\begin{split}
& \vert\varphi_\epsilon(x)-\varphi(x)\vert\leq\vert\varphi(x)\vert\leq\epsilon \ , \\ \\
& \vert\varphi_\epsilon^\prime(x)-\varphi^\prime(x)\vert\leq\vert\varphi^\prime(x)\vert
+(\delta_{x_\eta,\epsilon})^{-1}\chi^\prime 
\left((\delta_{x_\eta,\epsilon})^{-1}(x-x_\eta+2\delta_{x_\eta,\epsilon})\right)
\vert\varphi(x)\vert\leq C\epsilon \ ,
\end{split}
\end{equation}
where $C=1+\max_s\chi^\prime(s)$.
Consider now the case $a_{2,p}=0$. In this case we have $\bar\alpha=\min_{x\in\widetilde{E}\cup[a_{1,p},x_p]}$ 
and therefore $\varphi(\bar\alpha)=\varphi^\prime(\bar\alpha)=0$.
Hence 
\begin{equation} \label{eq620}
\varphi_\epsilon(x)=\varphi(x) \ , \ x\in[0,\bar\alpha) \ , 
\end{equation}
yields a well defined $C^1$ map and, again from \eqref{eq616} and \eqref{eq613}, we obtain
\begin{equation} \label{eq621}
\|(\varphi_\epsilon-\varphi)\vert_{[0,\bar\beta]}\|_{C^1[0,\bar\beta]}\leq C\epsilon \ . 
\end{equation}
By repeating the above construction for each connected component of  $\cup_{j=1}^NI_{2,j}$ we end up with a map
$\varphi_\epsilon\in C^1[0,1]$ that satisfies  $\|\varphi_\epsilon-\varphi\|_{C^1[0,1]}\leq C\epsilon$.
Since by construction $\varphi_\epsilon-\varphi\in C_{00}^1$, we have that $\varphi_\epsilon\in D(A)$.
Moreover, $\varphi_\epsilon$ has only a finite number of zeros and
\begin{equation} \label{eq622}
\varphi_\epsilon(x)\neq 0 \ \Rightarrow \ \varphi(x)\neq 0 \quad
\text{and}\quad\mathrm{sign}(\varphi_\epsilon(x))=\mathrm{sign}(\varphi(x)) \ .
\end{equation}
The completes the proof.
\end{proof}

\begin{lemma} \label{le63}
Assume the same as in Theorem \ref{th61}. Then
\begin{enumerate}
\item $\varphi(x)=\varphi^\prime(x)=0, \ \Rightarrow \ (A\varphi)(x)=0, \ \varphi\in D(A)$.
\item  For $\varphi\in D(A)$, $x\in(0,1)$, $(A\varphi)(x)$ depends only on the $1$-jet of $\varphi$ at $x$.
\item There exist $b,c\in C^0(0,1)$ such that
\begin{equation} \label{eq623}
(A\varphi)(x)=b(x)\varphi^\prime(x)+c(x)\varphi(x) \ , \ x\in(0,1) \ , \ \varphi\in D(A) \ .
\end{equation}
\end{enumerate}
\end{lemma}

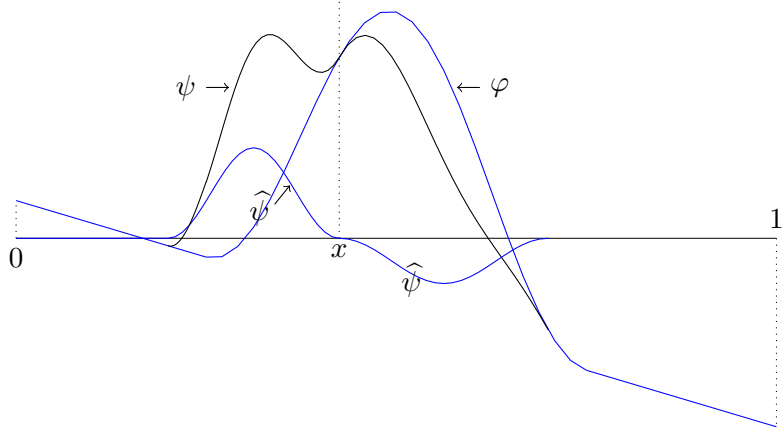
\begin{figure}
  \begin{center}
\begin{tikzpicture}[scale=1]
\draw[very thin] (-5,0)--(5,0);

\draw[blue, domain=-2.5:2.5] plot (\x, {2*(1-cos(72*(\x+2.5))-.5-.15*\x});
\draw[blue, domain=2.5:5] plot (\x, {2*(-.5-.15*\x)});
\draw[blue, domain=-5:-2.5] plot (\x, {2*(-.5-.15*\x)});
\draw[blue, domain=-3:-.75] plot (\x, {.6*(1-cos(160*(\x+3))});

\draw[ domain=-2.5:-.75] plot (\x, {2*(1-cos(72*(\x+2.5))-.5-.15*\x+.6*(1-cos(160*(\x+3))});
\draw[ domain=-.75:2] plot (\x, {2*(1-cos(72*(\x+2.5))-.5-.15*\x-.3*(1-cos(130.9*(\x+.75))});
\draw[ domain=-3:-2.5] plot (\x, {2*(-.5-.15*\x)+1*(1-cos(180*(\x+3))});



\draw[blue, domain=-.75:2] plot (\x, {-.3*(1-cos(130.9*(\x+.75))});

\draw[blue] (-5,0)--(-3,0);
\draw[dotted] (-.75,0)--(-.75,3.2);

\node[below] at (-.75,.05) {$x$};
\draw[dotted] (-5,0)--(-5,.5);
\draw[dotted] (5,0)--(5,-2.5);
\draw[->] (1.1,2)--(.8,2);
\draw[->] (-2.5,2)--(-2.2,2);
\node[right] at (1.1,2) {$\varphi$};
\node[left] at (-2.5,2) {$\psi$};
\draw[->] (-1.6,.45)--(-1.4,.7);
\node[] at (.2,-.55) {$\widehat{\psi}$};
\node[] at (-1.8,.40) {$\widehat{\psi}$};
\node[below] at  (-5,0) {$0$};
\node[above] at  (5,0) {$1$};
\end{tikzpicture}
\end{center}
\caption{The maps $\varphi$, $\widehat{\psi}$ and $\psi$.}
\label{fig3}
\end{figure}

\begin{proof}
To prove (i) we use an argument similar to the one used in the definition of the function $\varphi_\epsilon$
in the proof of Lemma \ref{le62}. For each $\delta>0$  sufficiently small we construct a map
$\varphi_\delta\in D(A)$ which vanishes in $(x-\delta,x+\delta)$ and satisfies
\begin{equation} \label{eq624}
\begin{split}
& \vert(\varphi-\varphi_\delta)(y)\vert\leq\vert\varphi(y)\vert\leq o(\delta) \ , \\ \\
& \vert(\varphi^\prime-\varphi_\delta^\prime)(y)\vert\leq\vert\varphi^\prime(y)\vert
+\delta^{-1}\chi^\prime\left(\delta^{-1}(y-x+2\delta)\right)\vert\varphi(y)\vert\leq f(\delta)
+\delta^{-1} o(\delta) \ .
\end{split}
\end{equation}

These estimates and the analogous in the intervals $[x-\delta,x+\delta]$ and $[x+\delta,x+2\delta]$ show that
$\varphi_\delta$ converges in $C^1[0,1]$ to $\varphi$ as $\delta\rightarrow 0$. This and assumption (iv) in
Theorem \ref{th61} imply $\lim_{\delta\rightarrow 0}\|A\varphi_\delta-A\varphi\|_{C^0[0,1]}=0$ and,
in particular, using also Lemma \ref{le62} which yields $(A\varphi_\delta)(x)=0$, we conclude $(A\varphi)(x)=0$.

To show (ii), let $x\in\supp(\varphi)$ be a point where $\varphi(x)$ and $\varphi^\prime(x)$ are not both zero.
Let $J_x^1\varphi\subset D(A)$ be the $1$-jet of $\varphi\in D(A)$ at $x\in(0,1)$.
From $\psi\in J_x^1\varphi$ it follows that $\widehat{\psi}=\psi-\varphi$ satisfies
$\widehat{\psi}(x)=\widehat{\psi}^\prime(x)=0$ (see Figure \ref{fig3}).
This and (i) imply $(A\widehat{\psi})(x)=0$ and therefore
\begin{equation} \label{eq625}
(A\psi)(x)=(A\varphi)(x)+(A\widehat{\psi})(x)=(A\varphi)(x) \ .
\end{equation}
This proves (ii).

From (ii) and the linearity of $A$ it follows that there is a continuous vector
$(0,1)\ni x\rightarrow \left((c(x),b(x)\right)$ such that
\begin{equation} \label{eq626} 
\begin{aligned}
(A\varphi)(x) &= \left((c(x),b(x)\right)\cdot\left(\begin{array}{l}\varphi(x) \\ 
\varphi^\prime(x)\end{array}\right) \\
&= c(x)\varphi(x)+b(x)\varphi^\prime(x) \quad , \quad x\in(0,1) \ , \;\varphi\in D(A) \ .
\end{aligned}
\end{equation}
The proof is complete.
\end{proof}

{\bf Remark 6.1:}
The functions $b$ and $c$ are defined and continuous in $(0,1)$. Since $A\varphi\in C^0[0,1]$ and
$(A\varphi)(x)=b(x)\varphi^\prime(x)+c(x)\varphi(x)$, $x\in(0,1)$, for all $\varphi\in D(A)$ it follows that
there exist
\begin{equation} \label{eq627}
\begin{split}
& L_0\varphi=\lim_{x\rightarrow 0^+}(A\varphi)(x)=b(0)\varphi^\prime(0)+c(0)\varphi(0) \ , \\ \\
& L_1\varphi=\lim_{x\rightarrow 1^-}(A\varphi)(x)=b(1)\varphi^\prime(1)+c(1)\varphi(1) \ ,
\end{split}
\end{equation}
and actually $b,c\in C^0[0,1]$.
For $\varphi\in D(A)$ and $0\leq t<<1$, $L_0$ and $L_1$ satisfy
\begin{equation} \label{eq628}
\begin{split}
& (T(t)\varphi-\varphi)(0)=tL_0\varphi+\emph{o}(t) \ , \\ \\
& (T(t)\varphi-\varphi)(1)=tL_1\varphi+\emph{o}(t) \ .
\end{split}
\end{equation}
The operators $L_0$ and $L_1$ may be defined for all $\varphi\in C^1[0,1]$  or impose some conditions
on $\varphi\in D(A)$ that restrict $D(A)$ to a proper subspace of $C^1[0,1]$.

From the previous analysis it follows that, if $T(t)$, $t\geq 0$, is a semigroup on $C^0[0,1]$ and the
domain $D(A)$ of the infinitesimal generator $A$ coincides with $C^1[0,1]$ (possibly restricted by the
conditions imposed by the operators $L_0$ and $L_1$), then there are functions $b,c\in C^0[0,1]$ such that
\begin{equation} \label{eq629}
(A\varphi)(x)=b(x)\varphi^\prime(x)+c(x)\varphi(x) \ , \ x\in[0,1] \ , \;\varphi\in C^1[0,1] \ .
\end{equation}
The definition of infinitesimal generator implies that the function $u$ defined by
$u(x,t)=(T(t)\varphi)(x)$ satisfies the scalar hyperbolic equation
\begin{equation} \label{eq630}
\left\{\begin{array}{l}
b(x)\frac{\partial}{\partial x}u(x,t)-\frac{\partial}{\partial t}u(x,t)=-c(x)u(x,t) \ , \\ \\
u(x,0)=\varphi(x) \ .
\end{array}\right.
\end{equation}
We assume $b(x)\neq 0$, $x\in[0,1]$ and solve this equation in case $b>0$. The case $b<0$ is similar.
The characteristic equations are
\begin{equation} \label{eq631}
\frac{d x}{d\tau}=b(x) \ , \ \frac{d t}{d\tau}=-1 \ , \ \frac{d \zeta}{d\tau}=-c(x)\zeta \ ,
\end{equation}
and it follows that 
\begin{equation} \label{eq632}
\int_{x_0}^x\frac{ds}{b(s)}+t=0 \ .
\end{equation}
Since $b>0$ by assumption, the integral in this equation is monotone increasing in $x$ and decreasing in $x_0$.
Furthermore, \eqref{eq632} can be solved for $x$ and $x_0$ defining functions $x=g(x_0,t)$ and $x_0=g^0(x,t)$ 
such that
\begin{equation} \label{eq633}
\int_{x_0}^{g(x_0,t)}\frac{ds}{b(s)}+t=\int_{g^0(x,t)}^x\frac{ds}{b(s)}+t=0 \ .
\end{equation}
We note the identity
\begin{equation} \label{eq634}
g(g^0(x,t),s)=g^0(x,t-s)=g(x_0,s) \ .
\end{equation}
This follows from \eqref{eq633} which implies
\begin{equation} \label{eq635}
\begin{split}
& \int_{g^0(x,t-s)}^x\frac{d\tau}{b(\tau)}+t-s=0 \ , \\
& \int_{g(g^0(x,t),s)}^x\frac{d\tau}{b(\tau)}+t-s=\int_{g^0(x,t)}^x\frac{d\tau}{b(\tau)}+t-
(\int_{g^0(x,t)}^{g(g^0(x,t),s)}\frac{d\tau}{b(\tau)}+s)=0 \ .
\end{split}
\end{equation}
The equation $x=g(x_0,t)$ is the equation of the characteristic or better of the projection of the
characteristic curve on the $x,t$ plane. By differentiating \eqref{eq633} with respect to $t$ one obtains
\begin{equation} \label{eq636}
\frac{\partial}{\partial t}g(x_0,t)=-b(g(x_0,t)) \ .
\end{equation}
This and $b>0$ show that the characteristic line oriented for increasing $t$ crosses the lines $x=Const$
from the right to the left. Moreover the characteristic line through $(x_0,t_0)$ is obtained from the
characteristic through $(x_0,0)$ by a translation of size $t_0$ in the time direction.
After setting $x=g(x_0,t)$ into the third equation \eqref{eq629} we obtain
\begin{equation} \label{eq637}
\zeta(x_0,t)=\zeta(x_0,0)e^{\int_0^tc(g(x_0,s))ds} \ .
\end{equation}
It follows that
\begin{equation} \label{eq638}
u(x,t)=\zeta(g^0(x,t),t)=\varphi(g^0(x,t))e^{\int_0^tc(g(g^0(x,t),s))ds}=
\varphi(g^0(x,t))e^{\int_0^tc(g^0(x,t-s))ds} \ ,
\end{equation}
where we have also used \eqref{eq634}. This equation defines the solution of \eqref{eq630} in the region
$R\subset[0,1]\times\R$ contained between the characteristic line through $(0,0)$  and the characteristic
line through $(1,0)$.
From \eqref{eq638} it follows that the sign of $u(x,t)$ along a characteristic curve $x=g(x_0,t)$ is constant
and coincides with the sign of $\varphi(x_0)$ where $x_0$ is the abscissa of the intersection of the considered
characteristic curve with the $x$ axes at time $t=0$. Indeed we have
\begin{equation} \label{eq639}
u(g(x_0,t),t)=\varphi(x_0)e^{\int_0^tc(g(x_0,s))ds} \ .
\end{equation}
Set $r=\int_0^1\frac{ds}{b(s)}$. Then, by means of \eqref{eq639}, we can define a linear homeomorphism
$h:C^0[-r,0]\rightarrow C^0[0,1]$ by setting (see Figure \ref{fig4})
\begin{equation} \label{eq640}
\begin{split}
&\psi:= u(1,\cdot) \ , \\
&h\psi\rightarrow\varphi \ .
\end{split}
\end{equation}
To see this, set $\xi(\theta)=g^0(1,\theta)$ and $\theta(\xi)=-\int_\xi^1\frac{ds}{b(s)}$.
The functions $\xi:[-r,0]\rightarrow[0,1]$ and $\theta:[0,1]\rightarrow[-r,0]$ are the inverse of each other:
\begin{equation} \label{eq641}
\xi(\theta(\xi))=\xi\in[0,1] \ ;\quad\;\theta(\xi(\theta))=\theta\in[-r,0] \ .
\end{equation}
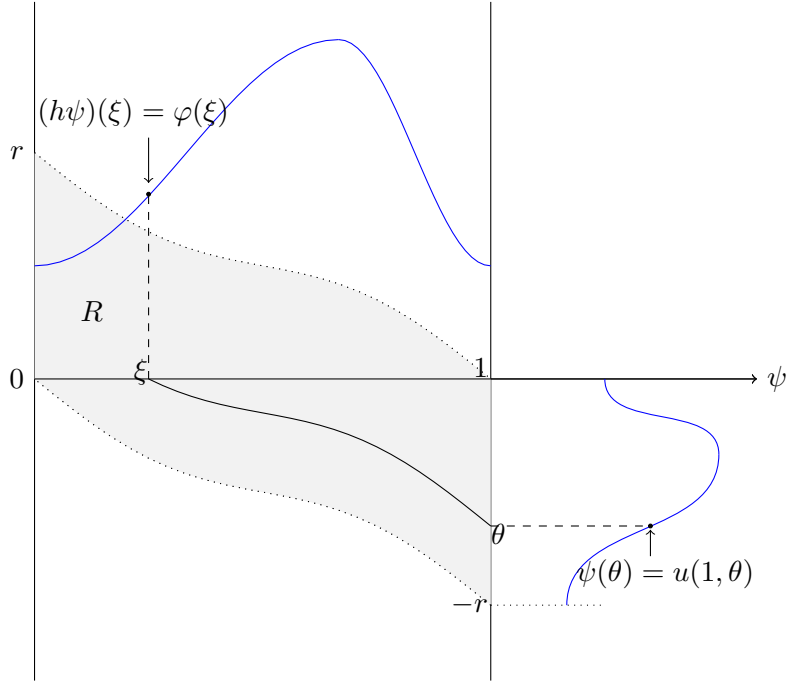
\begin{figure}
\begin{center}
\begin{tikzpicture}[scale=1]
\draw [->](0,0)--(9.5,0);
\draw [](0,-4)--(0,5);
\draw [](6,-4)--(6,5);
\draw [dotted](6,-3)--(7.5,-3);
\draw[dotted,domain=0:6] plot (\x,{.5*(6-\x)-.3*sin(60*\x)});
\draw[dotted,domain=6:0] plot (\x,{.5*(6-\x)-.3*sin(60*\x)-3});
\path [fill=lightgray!70!white!30!] (0,3)  [domain=0:6] plot (\x,{.5*(6-\x)-.3*sin(60*\x)}) --(0,0)--(6,0)--(6,-3)
[domain=6:0] plot (\x,{.5*(6-\x)-.3*sin(60*\x)-3})--(0,0)--(0,3);
\draw[dotted,domain=0:6] plot (\x,{.5*(6-\x)-.3*sin(60*\x)});
\draw[dotted,domain=6:0] plot (\x,{.5*(6-\x)-.3*sin(60*\x)-3});
\draw [->](0,0)--(9.5,0);
\draw[domain=1.5:6] plot (\x,{.5*(6-\x)-.3*sin(60*\x)-1.95});
\draw[blue] (7.5,0) to [out=270,in=90] (9,-1);
\draw[blue] (9,-1) to [out=270,in=90] (7,-3);
\draw[blue,domain=0:4] plot (\x,{1.5+1.5*(1-cos(45*\x))});
\draw[blue,domain=4:6] plot (\x,{4.5-1.5*(1-cos(90*\x))});
\draw [dashed](6,-1.95)--(8,-1.95);
\draw [dashed](1.5,0)--(1.5,2.5);
\node[right]at (9.5,0) {$\psi$};
\node[left]at (0,0) {$0$};
\node[left]at (6.1,.15) {$1$};
\node[left]at (0,3) {$r$};
\node[left]at (6.1,-3) {$-r$};
\node[]at (.75,.9) {$R$};
\node[]at (6.1,-2.05) {$\theta$};
\node[]at (1.4,.1) {$\xi$};
\draw[->] (8.1,-2.35)--(8.1,-2);
\draw[->] (1.5,3.2)--(1.5,2.6);
\node[above]at (1.3,3.2) {$(h\psi)(\xi)=\varphi(\xi)$};
\node[below]at  (8.3,-2.25) {$\psi(\theta)=u(1,\theta)$};
\draw [fill] (1.5,2.45) circle [radius=0.025];
\draw [fill] (8.1,-1.95) circle [radius=0.025];
\end{tikzpicture}
\end{center}
\caption{Illustration, for $b>0$, of $\psi(\theta)=u(1,\theta)$ and $\varphi(\xi)=h\psi(\xi)$. }
\label{fig4}
\end{figure}
Note also that \eqref{eq634} implies
\begin{equation} \label{eq642}
g(\xi(\theta),s)=g(g^0(1,\theta),s)=g^0(1,\theta-s)=\xi(\theta-s) \ .
\end{equation}
From \eqref{eq639} we have
\begin{equation} \label{eq643}
\begin{split}
& \psi(\theta(\xi))=u(1,\theta(\xi))=(h\psi)(\xi)e^{\int_0^{\theta(\xi)}c(g(\xi,s))ds} \\
& \Rightarrow \ (h\psi)(\xi)=\psi(\theta(\xi))e^{\int_{\theta(\xi)}^0c(g(\xi,s))ds} \ ,
\end{split}
\end{equation}
and
\begin{equation} \label{eq644}
(h^{-1}\varphi)(\theta)=u(1,\theta)=\varphi(\xi(\theta))e^{\int_0^\theta c(\xi(\theta-s))ds} \ ,
\end{equation}
where we have also used \eqref{eq642}.

\vskip.2cm

Equation \eqref{eq638} determines $(T(t)\varphi)(x)=u(x,t)$, $t\geq 0$ only in the set
$R\cap[0,1]\times[0,+\infty)$. To extend the solution to the whole of $[0,1]\times[0,+\infty)$ we need to
analyze the role and the structure of the boundary operators $L_0$ and $L_1$.
The fact that the solution is already defined on $R\cap[0,1]\times[0,+\infty)$ implies that $L_0$ is defined
for all $\varphi\in C^1[0,1]$. Instead $L_1:D(A)\rightarrow C^0[0,1]$ may impose effective conditions on $\varphi$
that restrict $D(A)$ to a proper subspace of $C^1[0,1]$.
In any case, if we set $y(t)=(T(t)\varphi)(1)$, then $y:[0,+\infty)\rightarrow\R$ is the solution of the problem
\begin{equation} \label{eq645}
\left\{\begin{array}{l}
y^\prime=L_1hy_t \ , \\ \\
y_0=h^{-1}\varphi \ ,
\end{array}\right.
\end{equation}
where
\begin{equation} \label{eq646}
y_t(\theta)=y(t+\theta) \ , \ \theta\in[-r,0] \ .
\end{equation}
Once $y$ is known, the semigroup is determined by
\begin{equation} \label{eq647}
T(t)\varphi=hy_t \ .
\end{equation}

{\bf Remark 6.2:}
We are under the assumption that $b>0$. If instead we assume $b<0$, then the operator $L_1$ does not impose
a restriction and can be defined for all $\varphi\in C^1[0,1]$ while the operator $L_0$ may restrict $D(A)$
to a proper subset of $C^1[0,1]$. Once the operator $L_0$ is specified, in case $b<0$, we can define the
analogous of the homeomorphism $h$ and show that the semigroup $T(t)$, $t\geq 0$ is determined by equations
similar to \eqref{eq645} and \eqref{eq647}.

We also observe that \eqref{eq645} includes delay equations as a special case. Indeed if we think of $T(t)$ 
as a $C_0$-semigroup on $C^0[-1,0]$ rather than on $C^0[0,1]$, then the delay equation corresponds to the 
case where $h$ is the identity.

\vskip.2cm

One can also consider the extension of \eqref{eq638} to negative time. For $b>0$ the extension to $t\leq 0$ 
depends on the choice of $L_0$ while, for $b<0$, depends on the choice of $L_1$.

If $b(x_0)=0$ for some $x_0\in(0,1)$, then $x=x_0$ is an asymptote for the characteristics. 
Moreover, $b(x)(x-x_0)>0$ for $x\neq x_0$ implies that the characteristics approach the asymptote 
for $t\rightarrow +\infty$. To extend the solution for $t\geq 0$ one needs to specify both 
operators $L_0$ and $L_1$ while the solution is automatically defined for $t\leq 0$. 
The opposite situation arises when  $b(x)(x-x_0)<0$, for $x\neq x_0$.
More complex situations where $b$ has several zeros can be analyzed along the lines developed before.
We discussed the case $b>0$ as the simplest significant case which, for $b\equiv 1$, includes 
the delay equation.

Next we determine what properties are required for $L_1$ in order that $V^\pm$ be a Lyapunov functional for
the semigroup defined by \eqref{eq647}.

\begin{lemma} \label{le64}
Let $L_1:D(A)\rightarrow\R$ be a  linear functional and assume that the semigroup $T(t)$, $t\geq 0$, defined by
\eqref{eq645} and \eqref{eq647} satisfies
\begin{equation} \label{eq648}
V^\mp(T(t)\varphi)\leq V^\mp(T(\tau)\varphi) \ , \quad t\geq\tau\geq 0 \ .
\end{equation}
Then
\begin{equation} \label{eq649}
L_1\varphi=a\varphi(0)+\alpha\varphi(1) \ ,
\end{equation}
with $\alpha\in\R$ and $\mp a\geq 0$, for $V^\mp$ respectively.
\end{lemma}

\begin{proof}
We only consider the case of the negative sign in \eqref{eq648}. The discussion of the other case is similar.
Observe that from \eqref{eq639} and the definition of $h$ it follows that 
\begin{equation} \label{eq650}
z(y_t)=z(hy_t) \ .
\end{equation}
Therefore it suffices to show that $V^-$ is a Lyapunov functional for the solution $t\rightarrow y_t$ of the
generalized delay equation \eqref{eq645}.
We divide the proof in tree steps:

\vskip.2cm

Step 1. Let $\varphi\in C^0[0,1]$ a map that satisfies $\varphi(1)=0$ and such that $z(\varphi)<+\infty$.
For small $t>0$, the solution $y_t$ of \eqref{eq645} satisfies
\begin{equation} \label{eq651}
y_t(\theta)=\left\{\begin{array}{l} (h^{-1}\varphi)(t+\theta) \ , \ \theta\in[-r+t,0] \ , \\ \\
\theta L_1\varphi+o(\theta) \ , \ \theta\in(0,t] \ ,
\end{array}\right.
\end{equation}
and, if $t>0$ is sufficiently small
\begin{equation} \label{eq652}
z(y_t\vert_{[-r+t,0]})=z(y_0)=z(\varphi) \ .
\end{equation}
The continuity of $\varphi$, \eqref{eq650} and the assumption $z(\varphi)<+\infty$ implies 
the existence of $\delta\in(0,r)$ and $\theta_0\in[-\delta,0]$ such that
\begin{equation} \label{eq653}
\mp y_0(\theta)\geq 0 \ , \ \theta\in[-\delta,0] \ \text{and} \ \mp y_0(\theta_0)>0 \ ,
\end{equation}
for $V^\mp$ respectively.
This and \eqref{eq652} imply that, if $\pm L_1(\varphi)>0$, then the function $y_t$ has a sign change in the
interval $[-\delta,t]$. Together with \eqref{eq408} this yields $z(y_t)=z(\varphi)+1$.
If $z(\varphi)$ is odd this implies $V^-(y_t)=V^-(y_0)+2$ in contradiction with \eqref{eq648}.
This contradiction proves that 
\begin{equation} \label{eq654}
\varphi\in C^0[0,1] \ , \ \varphi(1)=0 \ \text{and} \ z(\varphi) \ \text{odd}
\quad\Rightarrow\quad \varphi(\xi(\theta_0))L_1\varphi\geq 0 \ .
\end{equation}

\begin{figure}
\begin{center}
\begin{tikzpicture}[xscale=.7,yscale=1]
\draw[->](-6,.5)--(-6,2.5);
\draw [] (-6,1) -- (1.,1);
\draw [blue] (-6,1) -- (-5,1);
\draw [blue] (0,1) -- (1,1);

\draw[blue, domain=-5:-4] plot (\x, {1+sin(180*(\x+5))});
\draw[blue, domain=-1:0] plot (\x, {1+sin(180*(\x+1))});
\draw[blue, dotted, domain=-4:-3.2] plot (\x, {1-sin(180*(\x+4))});
\draw[blue, dotted, domain=-1.8:-1] plot (\x, {1+sin(180*(\x+1))});
\node[left] at (-6,1) {$0$};
\node[right] at (1,1) {$1$};
\node[left] at (-6,2) {$\varphi$};
\node[below] at (-5,1) {$\xi_1$};
\node[below] at (0,1) {$\xi_2$};
\draw[->](-6,-2.5)--(-6,-.5);
\draw [] (-6,-2) -- (1.,-2);
\draw [] (-6,-2) -- (1,-2);
\draw [blue]  (-5,-2) -- (0,-2);
\draw[blue, domain=-6:-5] plot (\x, {-2+1.5*sin(90-90*(\x+6))});
\draw[blue, domain=0:1] plot (\x, {-2-1.5*sin(180*(\x))});
\node[left] at (-6,-2) {$0$};
\node[right] at (1,-2) {$1$};
\node[left] at (-6,-1) {$\psi$};
\node[below] at (-5,-2) {$\xi_1$};
\node[below] at (0,-2) {$\xi_2$};
\draw [] (6,1) -- (13,1);
\draw[->](6,.5)--(6,2.5);
\draw [] (6,1) -- (13.,1);

\draw[blue, domain=7:8] plot (\x, {1+sin(180*(\x-7))});
\draw[blue, domain=11:12] plot (\x, {1+sin(180*(\x-11))});
\draw[blue, dotted, domain=8:8.8] plot (\x, {1-sin(180*(\x-8))});
\draw[blue, dotted, domain=10.2:11] plot (\x, {1+sin(180*(\x-11))});

\draw[blue, domain=6:7] plot (\x, {1+1.5*sin(90-90*(\x-6))});
\draw[blue, domain=12:13] plot (\x, {1-1.5*sin(180*(\x-12))});
\node[left] at (6,1) {$0$};
\node[right] at (13,1) {$1$};
\node[left] at (6,2) {$\varphi_\lambda$};
\draw [] (6,-2) -- (13,-2);
\draw[->](6,-2.5)--(6,-.5);
\draw [] (6,-2) -- (13,-2);
\draw[blue, domain=6:7] plot (\x, {-2+1.5*sin(90-90*(\x-6))});
\draw[blue, domain=12:13] plot (\x, {-2-1.5*sin(180*(\x-12))});

\draw[blue, domain=7:8] plot (\x, {-2-sin(180*(\x-7))});
\draw[blue, domain=11:12] plot (\x, {-2-sin(180*(\x-11))});
\draw[blue, dotted, domain=8:8.8] plot (\x, {-2+sin(180*(\x-8))});
\draw[blue, dotted, domain=10.2:11] plot (\x, {-2-sin(180*(\x-11))});
\node[left] at (6,-2) {$0$};
\node[right] at (13,-2) {$1$};
\node[left] at (6,-1) {$\varphi_\lambda$};
\end{tikzpicture}
\end{center}
\caption{$\varphi$, $\psi$ and $\varphi_\lambda$, for $\lambda>0$ and $\lambda<0$.}
\label{fig5}
\end{figure}
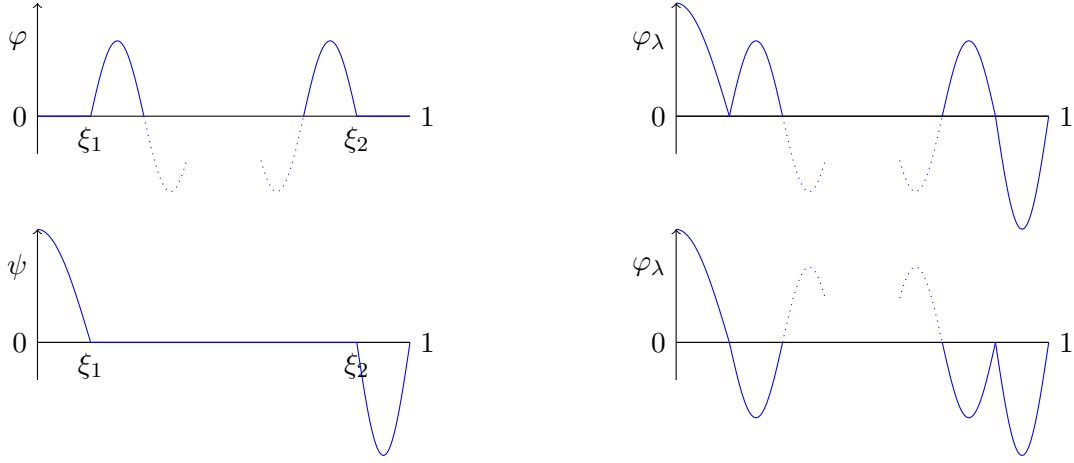

Step 2. Assume now that $\varphi\in C_{00}^0$ has $z(\varphi)<+\infty$ and even. Let  $0<\xi_1<\xi_2<1$
be such that $\supp(\varphi)\subset[\xi_1,\xi_2]$ and let $\psi\in C^0[0,1]$ be a function with
support in $[0,\xi_1]\cup[\xi_2,1]$ that satisfies $\psi(\xi)>0,\;\xi\in[0,\xi_1)$, 
$\psi(\xi)<0,\;\xi\in(\xi_2,1)$ and $\psi(1)=0$ (see Figure \ref{fig5}). The function $\psi$ has $z(\psi)=1$.
For $\lambda\in\R$ define $\varphi_\lambda=\psi+\lambda\varphi\in C^0[0,1]$. 
It results that $\varphi_\lambda(1)=0$ and
\begin{equation} \label{eq655}
z(\varphi_\lambda)=z(\psi)+z(\varphi)=z(\varphi)+1 \ , \ \lambda\in\R\setminus\{0\} \ .
\end{equation}
Hence $z(\varphi_\lambda)$ is odd for all $\lambda\in\R$ (see Figure \ref{fig4}).
Then \eqref{eq654} implies
\begin{equation} \label{eq656}
L_1\varphi_\lambda=\lambda L_1\varphi+L_1\psi\leq 0 \ , \ \lambda\in\R \ ,
\end{equation}
since $\varphi_\lambda(\xi(\theta_0))<0$.
If $L_1\varphi\neq 0$ then there exists $\lambda\in\R$ such that $L_1(\psi)+\lambda L_1\varphi>0$,
in contradiction with \eqref{eq654}, and therefore we conclude that
\begin{equation} \label{eq657}
\varphi\in C_{00}^0 \ \text{and} \ z(\varphi) \ \text{even} \ \Rightarrow \ L_1\varphi=0 \ .
\end{equation}
Since the set of functions with these properties is dense in the set
$C_0^0 = \{ \varphi\in C^0[0,1]: \varphi(0)=\varphi(1)=0 \}$, we have
\begin{equation} \label{eq658}
L_1\varphi=0 \ , \ \varphi\in C_0^0 \ ,
\end{equation}
whatever the parity of $z(\varphi)$.

\vskip.2cm

Step 3. Each $\varphi\in C^0[0,1]$  has a unique decomposition of the form
\begin{equation} \label{eq659}
\varphi=\widetilde\varphi+\varphi(0)\varphi^*+\varphi(1)\varphi_* \ , \ \widetilde\varphi\in C_0^0 \ ,
\end{equation}
where $\varphi_*(\xi)=\xi$, $\varphi^*(\xi)=1-\xi$, $\xi\in[0,1]$. This and \eqref{eq658} imply that 
\begin{equation} \label{eq660}
L_1\varphi=L_1\widetilde\varphi+\varphi(0)L_1\varphi^*+\varphi(1)L_1\varphi_*=a\varphi(0)+
\widetilde\alpha\varphi(1) \ , \ \varphi\in C^0[0,1] \ ,
\end{equation}
where we have set $a=L_1\varphi^*$ and $\widetilde\alpha=L_1\varphi_*$.
Finally, we observe that if $\varphi(1)=0$, $\varphi(0)\neq 0$ and $z(\varphi)$ is odd we have 
$\varphi(0)\varphi(\xi(\theta_0))<0$.
From this and \eqref{eq654} it follows that 
\begin{equation} \label{eq661}
\varphi(\xi(\theta_0))L_1\varphi\geq 0 \ \Rightarrow \ \varphi(0)L_1\varphi=\varphi(0)^2a\leq 0
\end{equation}
which implies $a\leq 0$ unless $L_1\varphi\equiv 0$. The proof is complete.
\end{proof}

\vskip.2cm

{\bf Remark 6.3:}
Note that the above analysis confirms the idea that relaxing the smoothness assumptions on the phase space 
restricts the class of semigroups which admit a given discrete Lyapunov functional. Indeed the discussion 
above shows that we must have $\alpha\equiv 0$ while $\alpha$ is allowed to be $\geq 0$ in Theorem 4.1.

\vskip.2cm

{\bf Remark 6.4:}
We can ask how the statement and the proof of Theorem \ref{th61} change when we replace (ii) by the assumption 
\begin{equation}\label{eq662}
z(T(t)\varphi)\leq z(T(\tau)\varphi) \ , \ \forall \ t\geq\tau \ , \ \varphi\in C^0[0,1] \ .
\end{equation} 
This replacement does not affect the proofs and the statements of Lemma \ref{le62} and Lemma \ref{le63}.
On the other hand, if \eqref{eq648} in Lemma \ref{le64} is replaced by \eqref{eq662}, then \eqref{eq649} becomes
\begin{equation}\label{eq663}
L_1\varphi=\alpha\varphi(1) \ , \ \alpha\in\R \ .
\end{equation} 
For the proof we observe that, with the new assumption, the argument in Step 1 yields 
\begin{equation}\label{eq664}
\varphi\in C^0[0,1] \ , \ \varphi(1)=0 \ \Rightarrow \ \varphi(\xi(\theta))L_1\varphi\geq 0 \ ,
\end{equation}
instead of \eqref{eq654}.
Then, arguing as in Step 2 but without accounting for the parity of $z(\varphi)$ one obtains \eqref{eq658} and 
\eqref{eq660}. Finally, if $\varphi(1)=0$, equations \eqref{eq664} and  \eqref{eq660} imply that 
$a\varphi(\xi(\theta))\varphi(0)\geq 0$ independently of $\varphi$. Therefore, $a=0$ follows.

\section{Conclusion/Discussion} \label{sec7}

The existence of a discrete Lyapunov functional has proved to be an essential tool for the detailed 
description of the global dynamics of equations \eqref{eq201} and \eqref{eq206}, 
\cite{firo96,firo99,firo00,fima89}. These results rise the question of existence of other classes 
of equations which admit some other type of discrete Lyapunov functionals.

A systematic study of this question would require the introduction of an abstract notion of discrete functional, 
say $\cV$ (see \cite{fulu97} for an attempt in this direction), which generalizes the known examples: the zero 
number $z$ and the related functionals $V^\pm$. Then, beginning with the finite dimensional case, one should 
determine all the pairs $(\cV,\cH)$ where $\cH$ is the class of linear operators which generates 
semigroups that admit $\cV$ as a discrete Lyapunov functional. This may be the object of future research but it 
is outside the scope of the present paper.
\vskip.1cm
For equations \eqref{eq201} and \eqref{eq206} the phase space is a set of scalar functions defined on an 
interval: a one dimensional domain. We guess that this fact is essential for the existence of a discrete 
Lyapunov functional. Therefore we conjecture that, besides the known examples, there is not too much to be 
discovered. In particular we don't expect the existence of a discrete Lyapunov functional for equations defined 
on higher dimensional domains or for vector valued equations like, for instance, the damped wave equation or 
parabolic equations on $n>1$ dimensional domains. See \cite{fulu97} where it is shown that no discrete 
functional exists for a parabolic equation on a two dimensional domain.

\vskip.2cm

In the search of other situations where the existence of a discrete Lyapunov functional can be ascertained 
we note that, if $H$ is a separable Hilbert space with the Hilbert basis $\{e_j\}_{j=1}^\infty$, the operators 
$A$ which in $\{e_j\}_{j=1}^\infty$ are represented by a semi-infinite Jacobi matrix with positive bounded off 
diagonal elements admit a discrete Lyapunov functional. But again this functional is related to the number of 
sign changes in the sequence of the coordinates of the generic element $x=\sum_{j=1}x_je_j$ of $H$. 
Note also that from Theorem 7.13 in \cite{sto32} any self-adjoint transformation with a simple spectrum has, 
in a suitable Hilbert basis, a matrix representation of semi-infinite Jacobi type with positive off diagonal 
elements. The study of this class of operators and their nonlinear counterparts may have a mathematical 
interest but does not seem to be related to some problem of physical relevance.

This seems to substantiate our conjecture that the discrete functionals defined via the zero number and 
their variations are actually the only possible. 

\vskip.2cm

It is well known that the parabolic equation \eqref{eq201} has a variational structure given by a \emph{continuous} 
Lyapunov functional, usually called \emph{energy function} \cite{zel68, mat88, firaro14,lafi19}. Therefore, a 
comment is in order to clarify the role played by  continuous and discrete Lyapunov functionals on the detailed 
description of the global dynamics defined by \eqref{eq201}. Under general dissipative conditions on the 
nonlinearity $f$, the decreasing character of the energy function leads to the existence of a global attractor 
$\cA$: a compact connected maximal set which is invariant under the dynamics defined by \eqref{eq201}. 
This is almost all the information that can be deduced from the existence of the energy function. 
On the other hand all the beautiful surprising results on the dynamics restricted to $\cA$, like the 
Morse-Smale property \cite{hen85,ang86,hal88} or the description of the global attractor $\cA$ from the 
\emph{ meander permutation} $\sigma\in S(n)$ \cite{furo91, firo99}, depend on the decay property of the zero 
number. This decay property induces, on the function space, a geometric structure consisting in a nested family 
of cones $K_n=\left\{ u\in C_n^1[0,1]:z(u)<n \right\}$, $n=1,\dots$. The key point is that the decay of the zero 
number implies the positively invariance of  $K_n$, $n=1,\dots$, both under the linearization of \eqref{eq201} 
around any of its solutions or for the difference of any two such solutions.

\vskip.2cm

In the previous Sections we have characterized the infinitesimal generators of semiflows with given
discrete Lyapunov functions derived from the zero number \eqref{eq202}. We have shown that, besides
the discrete Lyapunov functions associated to the semiflows, the smoothness required for the domain of
the infinitesimal generators also determines their characterization. In fact, we have shown that infinitesimal
generators acting on $C^2$ correspond to scalar parabolic equations while infinitesimal generators acting on
$C^1$ correspond to transport equations, both with adequate boundary conditions.

We next discuss some topics related to our results and some lines of research that could follow-up our research.

\vskip.2cm

Here we have discussed only equations with {\em discrete} Lyapunov functions. However, most equations 
considered here are also gradient-like, exhibiting a variational structure given by a 
{\em continuous} Lyapunov function. For example, this holds for scalar one-dimensional semilinear 
parabolic equations, with $C^2$-smooth dissipative nonlinearities, under separated boundary conditions, 
see \cite{zel68,mat88,firo14}. 
We remark that such equations define semiflows in $X=C^r$, $1<r<2$, \cite{olkr61,lasour68,lun95}. 
Note that due to the dissipative conditions on $f$, each semiflow has a global attractor $\cA$. 
Then, assuming hyperbolicity of all the equilibria, the {\em Sturm attractor} $\cA$ has the Morse-Smale property. 
See, for example, \cite{hen85,ang86,hal88}. This argument was also adapted to the case of periodic boundary 
conditions, see \cite{firowo04,firowo12a,firowo12b,firaro14}.
The existence of a continuous Lyapunov function also holds for quasilinear parabolic equations, 
see \cite{lap18}. In addition, we remark that Lappicy and Fiedler extended this result to the case of fully 
nonlinear scalar parabolic equations, see \cite{lafi19}. 

\vskip.2cm 

The continuous Lyapunov function is essential to establish the existence of the global attractor $\cA$, 
but from the mere knowledge of the continuous Lyapunov function very little information on the structure 
of $\cA$ could be derived. On the other hand, from the discrete Lyapunov function one could deduce which 
equilibria are connected by heteroclinic orbits, and the transversal intersection of stable and unstable 
manifolds of equilibria, see \cite{furo91,firo96,wol02}. 
We note that all these examples also exhibit discrete Lyapunov functions (lap numbers derived from the 
zero number). This also holds for the monotone feedback delay differential equations considered.

\vskip.2cm 

Finally, we remark that degenerate parabolic equations like $p$-Laplacian problems already mentioned 
in the Introduction, generate semiflows in $W^{1,p}(0,1)$, $p>2$, which possess a continuous Lyapunov function 
and a discrete Lyapunov function (i.e., the zero number), \cite{gebr07}.
This example raises the following questions: Can we obtain characterizations of infinitesimal generators acting
on less regular spaces $C^r[0,1]$ with $r<1$? Also, can we obtain such a characterization for fractional
semilinear parabolic equations? How far can we reduce $r$? We leave these questions to future research.

\section{Declaration} \label{sec8}

On behalf of the authors the corresponding author declares no conflict of interest.


\end{document}